\documentclass[a4paper,12pt]{article}
\usepackage{bbm}
\usepackage{amsmath}
\usepackage{amssymb}
\usepackage{amsfonts}
\usepackage{cite}
\usepackage[dvips]{graphicx}

\numberwithin{equation}{section}

\def \mff{\mathsf}

\newtheorem{theorem}{Theorem}[section]

\newtheorem{conjecture}[theorem]{Conjecture}
\newtheorem{corollary}[theorem]{Corollary}
\newtheorem{defn}[theorem]{Definition}

\newtheorem{lemma}[theorem]{Lemma}
\newtheorem{prop}[theorem]{Proposition}

\newtheorem{remark}[theorem]{Remark}

\def \df{\square}
\def \mc{\mathcal}
\def \inv{^{-1}}

\def \v{\vskip 0.1in}
\def \n{\noindent}

\def \real{\mathbb{R}}
\def \cplane{\mathbb{C}}
\def \indexn{^{(k)}}
\def \integer{\mathbb{Z}}
\def \detf{\det(f_{ s\bar{t}})}

\def \half{\frac{1}{2}}

\def \bb {\mathbb}
\def \b{\bar }
\def \bx{\mathrm{Box}}
\def \p{\partial}

\def \halfplane{\mathsf{h}^\ast}
\def \indexm{{_k}}

 \def\CHART{\mathsf{U}}

\def\t{\mathfrak{t}}

\title{Extremal metrics on toric surfaces}

\author{Bohui Chen, An-Min Li\footnote{Corresponding author.}, Li Sheng}

\begin{document}

\maketitle

{\abstract
In this paper, we study the Abreu equation on toric surfaces with prescribed  scalar curvatures on Delzant ploytope. In particular,  we
prove the existence of the extremal metric  when relative $K$-stability is
assumed.
\endabstract
}

\tableofcontents

\section{Introduction}

\v\n Extremal metrics, introduced by E. Calabi, have been
studied intensively in the past 20 years. There are three aspects of
the topic: sufficient conditions for existence, necessary
conditions for existence and the uniqueness of extremal metric.  The necessary conditions
for the existence are conjectured to be related to certain
stabilities. There are many works on this aspect
 (see \cite{CT,D1,D2,T,T1}). The uniqueness  was completed by
Mabuchi (see \cite{M})  and  Chen-Tian (cf.\cite{CT,CT1}).

On the other hand, there has been no much progress on the existence
of extremal metrics or K\"ahler metrics of constant scalar
curvature. One reason is that the equation is highly nonlinear and
of fourth order. The problem is to solve the equation under certain
necessary stability conditions. It was Tian who first  gave an analytic ``stability"
condition and  showed that such stability is equivalent to the existence of a K\"ahler-Einstein metric (see \cite{T1}).
In \cite{T1}, Tian also defined the algebro-geometric notion of $K$-stability. Then
in \cite{D2}, Donaldson generalized Tian's definition of $K$-stablilty by giving
an algebro-geometric definition of the Futaki invariant and conjectured that
it is equivalent to the existence of a  constant scalar curvature K\"ahler metric (cscK metric).
The problem may become  simpler if the manifold
admits more symmetry.
Hence, it is natural to consider the problem
on toric varieties first.
 Since each toric manifold $M^{2n}$ can be represented by a Delzant
polytope in $\real^n$, by Abreu, Burns and Guillemin's work, the fourth order equation can be transformed to be an equation of
real convex functions on the polytope, which is known as the
Abreu equation.
In a sequence of papers, Donaldson
initiated  a program to study the extremal metrics on toric manifolds:   Donaldson
formulated  $K$-stability for polytopes and conjectured  that the
stability implies the existence of the cscK metric on toric manifolds.
 In \cite{D2}, Donaldson also proposed a
stronger version of stability which we call   {\em uniform
stability} in this paper (cf.\;Definition \ref{definition_2.1.5}).
 The existence of weak solutions was solved by Donaldson under the assumption of the
uniform   stability (cf. \cite{D2}). Note that in \cite{CLS5}, we prove that the uniform
stability is a necessary condition. On the other hand,
Zhou-Zhu (cf.\cite{ZZ1}) introduced the notion of properness on the modified
Mabuchi functional and showed the existence of weak solutions under
this assumption. When the scalar curvature $K>0$ and $n=2$, all these
conditions are equivalent. The remaining issue is to show the
regularity of the weak solutions. In a sequence of papers (cf.
\cite{D2,D3,D4}), Donaldson solved the problem for cscK metrics on toric surfaces by proving the
regularity of the weak solutions.

 In this paper, we study the existence for metrics of {\em any
prescribed scalar curvature} on Delzant polytope (including extremal metrics) under the assumption
of uniform   stability.  As we know, though the equation of extremal metric is on the complex
manifolds, for the toric manifolds, the equation can be reduced to a
real equation on the Delzant polytope. The second
author with his collaborators developed a systematic package of tools to study one type
of fourth order PDEs which includes the Abreu
equation (cf.\cite{L-J,L-J-1,L-J-2,L-J-3,L-J-4,L-X1,L-X2,L-X-S-J,CLS1,CLS2,L-S}).
We call the package  the {\em real affine techinique}. This is explained
in \cite{CLS}. The challenging part   is then to study the
boundary behavior of the Abreu equation near the boundary of
polytope. The boundary of polytopes can be thought as the
interior of the complex manifold.   The important  issue is then to
generalize the affine techniques to the complex case. In \cite{CLS},
we made an attempt on this direction.
In particular,
we obtain the estimate of the Ricci curvature $\mc K$  in terms of the bound of $H$ (cf. Theorem \ref{theorem_3.4.2}).
We call the technique  the {\em complex affine technique}.
The real/complex affine techniques play an important role
in both  \cite{CLS} and  this paper.

 The main theorem of this paper is
\begin{theorem}\label{theorem_1.1}
Let $(M,\omega)$ be a compact toric surface and $\Delta$ be its
Delzant polytope. Let $K\in C^{\infty}(\bar\Delta)$ be an edge-nonvanishing
function. If $(\Delta,K)$ is uniformly stable, then there is a smooth
$T^2$-invariant metric $\omega_{f}$ on $M$ such that the scalar curvature of $\omega_{f}$ is $K$.
\end{theorem}
For the definitions of the edge-nonvanishing function and  uniform
stability, readers are  referred to Definition \ref{definition_1} and
Definition \ref{definition_2.1.5} respectively. For any Delzant
polytope $\Delta,$ there is a unique affine linear function $K$ such that
$\mathcal L_K(u)=0$  for all affine linear functions $u$ (for the
definition of $\mathcal L_K $ see \S3.1). On the other hand, by a
result of Donaldson (cf. Proposition \ref{proposition_2.1.6}), we
know that any  relative $K$-polystable $(\Delta,K)$ with $K>0$ is
uniformly stable.  Hence, as a corollary, we get the following
theorem.
\begin{theorem}\label{theorem_1.2}
Let $M$ be a compact toric surface and $\Delta$ be its Delzant
polytope. Let $K\in C^{\infty}(\bar\Delta)$ be a positive  linear function. If $(M,K)$ is relative K-polystable then there is a smooth
$T^2$-invariant metric on $M$ with scalar curvature  $K$.
\end{theorem}

The paper is organized as follows: in \S \ref{sect_1}, we review the background and formulate the problems.
In \S \ref{sect_2.2}, we explain that Theorem
\ref{theorem_1.1} can be reduced to Theorem \ref{theorem_2.2.2}.
The proof of
Theorem \ref{theorem_2.2.2} occupies the rest of the paper and the contents of sections are indicated by the titles.

Acknowledgments. We would like to thank Xiuxiong Chen, Yongbin Ruan and Gang Tian
for their interests in the project and continuing supports. We would also like to thank
Bo Guan  and
Qing Han  for their  valuable discussions. The first author is
 partially supported
by NSFC 10631050 and NKBRPC (2006CB805905), and the second author is
partially supported by NSFC 10631050, NKBRPC (2006CB805905) and RFDP (20060610004).

%%%%%%%%%%%%%sect_1%%%%%%%%%%%%%%%%%%

\section{K\"{a}hler geometry on toric surfaces }\label{sect_1}

%%%%%%%%%%%%%%%%%%%%%%%%sect_2

 In this section, we review  the K\"{a}hler geometry {of} toric surface and introduce the notations to be used in this paper. We assume that the readers are familiar with some basic knowledge of toric varieties.

A toric manifold is a 2$n$-dimensional K\"{a}hler manifold
$(M,\omega)$ that admits an $n$-torus (denoted $\bb T^n$) Hamitonian action.
Let $\tau: M\to \mathfrak{t}^\ast$ be the moment map, where
$\t^\ast$, identified with $  \real^n$, is the  dual of the  $\mathfrak{t}$ which is the Lie algebra of $\bb T^n$. The image $\bar\Delta=\tau(M)$ is known to be
a polytope(\cite{De}).
In the literature, people use $\Delta$ for the image of the moment map.
However, it is  more convenient in this paper to always assume that
{\em $\Delta$ is an open polytope.}  Note that
 $\Delta$ determines a fan $\Sigma$ in $\t$.
The converse is not true: $\Sigma$ determines $\Delta$ only up to a
certain similarity. $M$ can be reconstructed from either $\Delta$ or
$\Sigma$(cf. \cite{Fu} and \cite{Gu}). Moreover, the
class of $\omega$ can be read from $\Delta$. Hence, the uncertainty
of $\Delta$ reflects the non-uniqueness of K\"{a}hler classes. Two
different constructions are related via Legendre transformations.
The K\"{a}hler geometry appears naturally when considering the
transformation between two different constructions. This was
explored by Guillemin in \cite{Gu}. We will summarize these facts in
this section. For simplicity, we only consider the toric surfaces,
i.e, $n=2$.

\subsection{Toric surfaces and coordinate charts}\label{sect_1.1}

Let $\Sigma$ and $\Delta$ be a pair consisting of a fan and a polytope for a
toric surface $M$. For simplicity, we focus on compact toric
surfaces. Then $\Delta$ is a Delzant polytope in $\t^\ast$ and
its closure is compact.

We use the notations in \S 2.5  of \cite{Fu} to describe the fan. Let
$\Sigma$ be a fan given by a sequence of lattice points
$$\{
\nu_0,\nu_1,\ldots, \nu_{d-1},\nu_d=\nu_0 \}$$ in  the counterclockwise order, in
$N=\integer^2\subset \t$, such that successive pairs generate
$N$.
Suppose that the vertices and  edges of $\Delta$ are denoted by
$$
\{\vartheta_0,\ldots,\vartheta_d=\vartheta_0\}, \;\;\;
\{\ell_0,\ell_1,\ldots, \ell_{d-1},\ell_d=\ell_0\}.
$$
Here $\vartheta_i=\ell_i\cap \ell_{i+1}$.

By saying that $\Sigma$ is dual to $\Delta$ we mean that $\nu_i$ is
the inward pointing normal vector to $\ell_i$ of $\Delta$. Hence,
$\Sigma$ is determined by $\Delta$. Suppose that the equation for
$\ell_i$ is
\begin{equation}\label{eqn_1.1}
l_i(\xi):=\langle\xi,\nu_i\rangle- \lambda_i=0.
\end{equation}
 Then  we have
$$
\Delta=\{\xi|l_i(\xi)> 0,\;\;\; 0\leq
i\leq d-1\}
$$

There are three types of cones in $\Sigma$:  a 0-dimensional cone
$\{0\}$  denoted by $C_\Delta$;   1-dimensional cones
generated by $\nu_i$ and denoted by  $C_{\ell_i}$;
2-dimensional cones generated by $\{\nu_{i},\nu_{i+1}\}$ and
denoted by $C_{\vartheta_i}$.
It is known that for each cone of $\Sigma$, one can associate  to it a
complex coordinate chart  of $M$ (cf. \S1.3 and \S 1.4 in
\cite{Fu}).
 Let $\CHART_\Delta,
\CHART_{\ell_i}$ and $\CHART_{\vartheta_i}$ be the coordinate
charts. Then
$$
\CHART_{\Delta}\cong (\cplane^\ast)^2;\;\;\;
\CHART_{{\ell_i}}\cong \cplane\times\cplane^\ast ;\;\;\;
\CHART_{{\vartheta_i}}\cong \cplane^2.
$$
In particular, in each $\CHART_{\ell_i}$ there is a divisor
$\{0\}\times \cplane^\ast$. Its closure is a divisor in $M$, we
denote it by $Z_{\ell_i}$.

\begin{remark}\label{remark_1.1.1}
$\cplane^\ast$ is called a complex torus and denoted by $\bb T^c$. Let $z$ be its natural
coordinate.

In this paper, we introduce another complex coordinate by
considering the following identification
\begin{equation}\label{eqn_1.1a}
 \bb T^c\to \real\times 2\sqrt{-1}\bb T;\;\;\;
 w=\log z^2.
\end{equation} We call $w=x+\sqrt{-1}y$  the
 log-affine complex coordinate (or  log-affine coordinate) of $\cplane^\ast$.

When $n=2$, we have
$$
(\cplane^\ast)^2\cong \t\times 2\sqrt{-1}\bb T^2.
$$
Then $(z_1,z_2)$ on the left hand side is  the usual complex
coordinate; while $(w_1,w_2)$ on the right hand side is  the
 log-affine coordinate.
Write $w_i=x_i+\sqrt{-1}y_i,$ $y_{i}\in [0,4\pi]$.
Then $(x_1,x_2)$ is the coordinate of $\t$.
\end{remark}

We make the following convention.
\begin{remark}\label{remark_1.1.2}
On different types of coordinate chart  we use different
coordinate systems as follows:
\begin{itemize}
\item on $\CHART_\vartheta\cong \cplane^2$, we use the coordinate
$(z_1,z_2)$; \item on $\CHART_\ell\cong
\cplane\times\cplane^\ast$, we use the coordinate $(z_1,w_2)$;
\item on $\CHART_\Delta\cong (\cplane^\ast)^2$, we use the
coordinate $(w_1,w_2)$, \end{itemize}
 where $z_i=e^{\frac{w_i}{2}},i=1,2$.
\end{remark}

\begin{remark}\label{remark_1.1.3}
Since we study the $\bb T^2$-invariant geometry on $M$, it is
useful to { specify} a representative point of each $\bb T^2$-orbit.
Hence for $(\cplane^\ast)^2$, we let the points on $\t\times
2\sqrt{-1}\{1\}$  be the representative points.
\end{remark}

\subsection{K\"{a}hler geometry on toric surfaces}\label{sect_1.2}
  Guillemin  in \cite{Gu} constructed  a natural $T^2$-invariant K\"{a}hler  form
$\omega_o$ on  $M$ from the polytope $\Delta$.  We take this form  as a reference point in  the class
 $[\omega_o]$ and call  the associated K\"{a}hler metric  the Guillemin metric.

For each $T^2$-invariant K\"{a}hler form $\omega\in [\omega_o]$, on each coordinate chart associated to a
cone of the fan $\Sigma$, there
is a   K\"{a}hler  potential function (unique up to
 linear functions).  Write the collection of the potential functions as
$$
\mathsf f=\{f_\bullet\}:=\{f_\Delta, f_{\ell_0},\ldots, f_{\ell_{d-1}}, f_{\vartheta_0},\ldots,
f_{\vartheta_{d-1}} \}
$$
as  the   K\"{a}hler potential function with respect to the coordinate charts
$$\{
\CHART_\Delta,
\CHART_{\ell_i}, \CHART_{\vartheta_i} | i=1, \cdots, d-1\}.$$
We write $\omega=\omega_{\mathsf f}$   to indicate the associated potential functions.

 Let
$ \mathsf{g}=\{g_\bullet \}$
be the collection of potential functions for $\omega_o$.  We can  realize $\mathsf f$ by the following construction.
 %Readers are referred  to the end of \S\ref{sect_1.3} for the construction of $\mathsf g$.
Let $C^\infty_{\bb T^2}(M)$ be the smooth $\bb T^2$-invariant
functions of $M$.  Set
$$
 C^\infty_+(M) =
\{\phi\in C^\infty_{\bb T^2}(M)|\omega_{g_\Delta+\phi}>0\}.$$
Then for $\omega\in [\omega_o]$ there exists a function $\phi$ such that
$$
\mathsf f=\mathsf g+\phi:=(g_\bullet+\phi)
$$
and $\omega=\omega_\mathsf f$. Set
$$
\mc C^\infty(M,\omega_o)=\{\mathsf f|\mathsf f=\mathsf g+\phi, \phi\in C^\infty_+(M)\}.
$$

\begin{remark}\label{remark_1.2.1}
 Suppose that  $f_\bullet=g_\bullet+\phi.$
Consider the matrix
$$
\mathfrak{ M}_ {f}=(\sum_k { g}^{i\bar
k} {f}_{j\bar k}).
$$
Though this is not a globally well defined matrix on $M$, its
eigenvalues are globally defined. Set $ \nu_f$ to be the set of eigenvalues and $
H_ {f}=\det\mathfrak{ M}_ {f}\inv. $ These are global
functions on $M$.
\end{remark}
 Under   a coordinate chart with potential function $f$, the
Christoffel symbols, the curvature tensors, the Ricci curvature
and the scalar curvature of K\"{a}hler metric
$\omega_{{f}}$ are given by
$$ {\Gamma}^{k}_{ij}=\sum_{l=1}^{n} f^{k \b l}\frac{\p f_{i \b l}}{\p z_j}, \;\;\;
 {\Gamma}^{\b k}_{\b i\b j}=\sum_{l=1}^{n} f^{\b k  l}\frac{\p
f_{\b i l}}{\p z_{\b j}} ,  $$
$$ {R}_{i\b jk\b l}=- \frac{\p ^2 f_{i\b j}}{\p z_k \p z_{\b l}}+\sum f^{p\b q}\frac{\p f_{i \b q}}{\p z_k}
\frac{\p f_{p \b j}}{\p z_{\b l}},$$
$$
R_{i\b j}= - \frac{\partial^{2}}{\partial z_{i}\partial \bar
z_{j}} \left(\log \det\left(f_{k\bar l}\right)\right), \;\;\;\mc
S=\sum f^{i\bar j}R_{i\b j},
$$
respectively. When we use the log-affine coordinates
on $\t$, the Ricci curvature
and the scalar curvature can be written as
$$
R_{i\b j}= - \frac{\partial^{2}}{\partial x_{i}\partial x_{j}}
\left(\log \det\left(f_{kl}\right)\right), \;\;\;\mc S=-\sum
f^{ij}\frac{\partial^{2}}{\partial x_{i}\partial x_{j}} \left(\log
\det\left(f_{ij}\right)\right).
$$
 We  treat
$\mc S$ as an operator for $f$ and denote it by
$\mc S(f)$.

 Define
\begin{equation}\label{eqn_1.1c}
\mc K\;=\;\|Ric\|_f +\|\nabla Ric\|_f
^{\frac{2}{3}}+\|\nabla^2 Ric\|_f
^{\frac{1}{2}}
\end{equation}
\begin{equation}
W=\det(f_{s\bar t}),\;\;\;\Psi=\|\nabla\log W\|^2_f.
\end{equation}
 We also  denote by $\dot \Gamma_{ij}^k$, $\dot
R^m_{ki\bar{l}}$  and $\dot R_{i\b j}$ the connections, the
curvatures and the Ricci curvature of the metric $\omega_{o} $
respectively.

When  focusing on $\CHART_\Delta$ and using the log-affine
coordinate (cf.\;Remark \ref{remark_1.1.1}), we have $f(x)=g(x)+\phi(x)$. We  remark
that when restricting on $\real^2\cong \real^2\times2\sqrt{-1}\{1\}$,
the Riemannian metric induced from $\omega_\mathsf{f}$ is
 the Calabi metric $G_f$
(cf. \S \ref{sect_2.1}).

Fix a large constant $K_o>0$.
 We set
$$
\mc C^\infty(M,\omega_o;K_o)
=\{\mathsf{f}\in \mc C^\infty(M,\omega_o)|
|\mc S(f)|\leq K_o\}.
$$
In this paper, we  mainly study the apriori estimates for the functions in this class.

\subsection{The Legendre transformation, moment maps and potential functions}
\label{sect_1.3}
 Let $f$ be a (smooth) strictly convex function on
$\t$.
 The gradient of $f$
defines a (normal) map $\nabla^f$ from $\t$ to  $\t^\ast$:
$$
\xi=(\xi_1,\xi_2)=\nabla^f(x) =\left(\frac{\partial f}{\partial
x_1},\frac{\partial f}{\partial x_2}\right).
$$
The function $u$ on $\t^\ast$
$$
u(\xi)=x\cdot\xi - f(x)
$$
is called the Legendre transformation of $f$. We write $u=L(f)$.
Conversely, $f=L(u)$.

Now we restrict  to $\CHART_\Delta.$ When we use the coordinate $(z_1,z_2),$
the moment map with respect to $\omega_{f}$
 is given by
\begin{equation}\label{eqn_1.1b}
\tau_{{\mathsf{f}}}:\CHART_\Delta\xrightarrow{(\log |z_1^2|,\log|z_2^2|)}
\t\xrightarrow{\nabla^f} \Delta
\end{equation}
Note that the first map is induced from \eqref{eqn_1.1a}.  It
is known that $u$ must satisfy certain behavior near boundary of
$\Delta$.
\begin{theorem}[Guillemin]\label{theorem_1.3.1}
{ Let $v=L(g)$, where $g=g_\Delta$ is the potential
function of the Guillemin metric. }Then
$v(\xi)=\sum_{i} l_i\log l_i$, where $l_i$ is defined in \eqref{eqn_1.1}.
\end{theorem}
	For $u=L(f)$,   we have $u=v+\psi$, where $\psi\in C^\infty(\bar\Delta)$. Motivated by this,
	  we set
$$
\mc C^\infty(\Delta,v)=\{u| u=v+\psi \mbox{ is strictly convex, }
\psi\in C^\infty(\bar \Delta)\},
$$
with $v$ as a
reference point.   Note that  this space only depends on $\Delta$.

We summarize the fact we just presented: let $\mathsf{f}\in \mc
C^\infty(M,\omega_o)$, then the moment map $\tau_{{\mathsf{f}}}$
is given by $f=f_\Delta$ via the diagram \eqref{eqn_1.1b} and $u=L(f)\in \mc
C^\infty(\Delta,v)$. Conversely, $f_\bullet$ can be constructed from $u$ as well.

Given a function $u\in \mc C^\infty(\Delta,v)$, we can get an
$\mathsf{f}\in \mc C^\infty(M,\omega_o)$ as { follows}.
\begin{itemize}
\item On $\CHART_\Delta$,
 $f_\Delta=L(u)$;
\item  on $\CHART_{\vartheta_i}$, $f_{\vartheta_i}$ is constructed in the
following steps: (i), since $\vartheta_i=\ell_i\cap \ell_{i+1}$, let $B\in SL(2,\mathbb Z)$ be the transformation of $\mathfrak t^\ast$ such that
$$
B(\nu_i)=(1,0), \;\;\;B(\nu_{i+1})=(0,1).
$$
 Meanwhile, $u$ is transformed to a function in the following format
$$u'=\xi_1\log\xi_1+\xi_2\log\xi_2+\psi';$$
(ii), $f'=L(u')$ defines  a function on $\t$ and therefore is a
function on $(\cplane^\ast)^2\subset
 \CHART_{\vartheta_i}$ in terms of log-affine coordinate; \\
 (iii),  it is known that
 $f'$ can be extended over $\CHART_{\vartheta_i}$ and we set $f_{\vartheta_i}$
to be this function; \item on $\CHART_\ell$, the construction of
$f_\ell$ is similar to $f_\vartheta$. The reader may refer to
\S\ref{sect_1.5} for the construction.
\end{itemize}

\subsection{The Abreu equation on $\Delta$}\label{sect_1.4}
We can transform  the scalar curvature operator $\mc S(f)$ to  an
operator $\mc S(u)$ of $u$ on $\Delta$, where $u=L(f)$.   Then
$$
\mc S(u)=\mc S(f)\circ\nabla^u.
$$
  The operator $\mc S(u)$  is known
to be
$$
\mc S(u)=-\sum U^{ij}\partial^{2}_{ij} w
$$
where $(U^{ij})$ is the cofactor matrix of the Hessian matrix
$(\partial^2_{ij}u)$, $w=(\det(\partial^2_{ij}u))\inv$.
Here and later we denote $\partial^2_{ij}u=\frac{\partial^2u}{\partial \xi_i\partial \xi_j}.$ It is well known that $\omega_{\mathsf{f}}$ gives an extremal
metric if and only if $\mc S\circ\nabla^ u$ is a linear function
of $\Delta$. Let $K$ be a smooth function on $\bar\Delta$, the Abreu
equation is
\begin{equation}\label{eqn_1.2}
\mc S(u)=K.
\end{equation}
We set $\mc C^\infty(\Delta, v;K_o)$
to be the functions {\bf $u\in\mc C^\infty(\Delta,v)$} with $|\mc S(u)|\leq K_o$.

\begin{defn}\label{definition_1}
Let $K$ be a smooth function on $\bar\Delta$. It is called {\em edge-nonvanishing}
if it does not vanish on any edge of $\Delta$. That is to say,
for any edge $\ell$ there exists a point $\xi^{(\ell)}$ on the edge such that
$K(\xi^{(\ell)})\not=0$.
\end{defn}
In our papers, we will always assume that $K$ is edge-nonvanishing.

\subsection{A special case: $\cplane\times \cplane^\ast$}\label{sect_1.5}

Let $\halfplane\subset \t^\ast$ be the half plane given by
$\xi_1\geq 0$. The boundary is { the} $\xi_2$-axis and we  denote it by
$\t_2^\ast$.
 The corresponding fan consists of only one lattice $\nu=(1,0)$.
 The coordinate chart is
$\CHART_{\halfplane}= \cplane\times\cplane^\ast$. Let
$Z=Z_{\t_2^\ast}=\{0\} \times \cplane^\ast$ be its divisor.

Let $v_{\halfplane}= \xi_1\log\xi_1+\xi_2^2$. Set
$$
\mc C^\infty(\halfplane,v_{\halfplane})=\{u|u=v_{\halfplane}+\psi
\mbox{ is strictly convex, }\psi\in C^\infty(\halfplane)\}
$$
and $\mc C^\infty(\halfplane, v_{\halfplane};K_o)$
be the functions whose $\mc S$ is less than $K_o$.

Take  a function $u\in \mc C^\infty(\halfplane,
v_{\halfplane})$. Then $f=L(u)$ is a function on $\t$.
Hence it defines a function on the $\cplane^\ast\times\cplane^\ast
\subset \CHART_{\halfplane}$ in terms of log-affine {  coordinates}
$(w_1,w_2)$. Then the function $  f_\mathsf{h}(z_1,w_2):=f(\log |z_1^2|,
Re(w_2)) $ extends smoothly over $Z$, {  and hence is defined} on
$\CHART_{\halfplane}$. We conclude that for any $u\in \mc
C^\infty(\halfplane,v_{\halfplane})$ it yields a potential
function $ f_\mathsf{h}$ on $\CHART_{\halfplane}$.

 When we choose the coordinate $(z_1,w_2),$
the moment map with respect to ${\omega_{ {f}}}$
 is given by
\begin{equation}\label{eqnc_1.1b}
\tau_{{\mathsf{f}}}:\CHART_\halfplane\xrightarrow{(\log |z_1^2|,Re(w_2))}
\t\xrightarrow{\nabla^f}\halfplane.
\end{equation}

Using $v_{\halfplane}$ and the above argument, we define a
function $g_\mathsf{h}$ on $\CHART_\halfplane$.

\subsection{$K$-stability }\label{sect_2.1}
%\subsection{$K$-stability and existence of extremal metrics}\label{sect_2}
In a sequence of papers, Donaldson initiates a program to study
the extremal metrics on toric manifolds. Here, we outline his
program and
 some
of his important results. Again, we restrict ourself only on the
2-dimensional case.

Let $\Delta$ be a Delzant
 polytope in $\t^*$.  Most of the material in this subsection can be applied
 to general convex polytopes, or even convex domains. However, for simplicity
 we focus on the Delzant polytopes.

For any smooth function $K$ on $\bar\Delta$, Donaldson defined  a
functional on $\mc C_\infty(\Delta)$:
$$\mc F_{K}(u)= -\int_{\Delta}\log \det(\partial^2_{ij}u)d\mu + \mc L_{K}(u),$$
where $\mc L_{K}$ is the linear functional
$$\mc L_{K}(u) = \int_{\partial \Delta}u d\sigma - \int_{\Delta}Ku
d\mu,$$
where $d\mu$ is the  Lebesgue measure on $\mathbb R^n$ and on each face F $d\sigma$ is a constant multiple of the standard $(n-1)$-dimensional Lebesgue measure (see \cite{D1} for details).
  In \cite{D1}, Donaldson defined the concept of
$K$-stability by using the test configuration (Definition 2.1.2 in
\cite{D1}). We recall the definition of relatively  $K$-polystability for
toric manifolds(cf. \cite{Sz}).
\begin{defn}\label{definition_2.1.1}{\it [relatively $K$-polystable]}
 Let $K\in
C^{\infty}(\bar\Delta)$ be a smooth function on $\bar\Delta$.
$({\Delta},K)$ is called {\it relatively $K$-polystable} if $\mc
L_{K}(u)\geq 0$ for all rational piecewise-linear convex functions
$u$,
 and $\mc
L_{K}(u)=0$ if and only if $u$ is a linear function.
\end{defn}
In this paper, we will simply refer to  relatively K-polystable as
polystable.

We fix a point $p\in \Delta$ and say $u$ is normalized at $p$ if
$$u(p)\geq 0,\;\; \nabla u(p)=0.$$
By Donaldson's work, we make the following definition.
\begin{defn}\label{definition_2.1.5}
$({\Delta},K)$ is called uniformly stable if it is polystable and
for any normalized convex function $u\in \mc C^\infty(\Delta)$
$$
\mc L_K(u)\geq \lambda\int_{\partial \Delta} u d \sigma
$$
for some constant $\lambda>0$. Sometimes, we say that $\Delta$ is
$(K,\lambda)$-stable.
\end{defn}
  Donaldson proved
\begin{prop}\label{proposition_2.1.6}
When $n=2$, if $(\Delta,K)$ is  polystable and $K>0$, then there exists a
constant $\lambda>0$ such that $\Delta$ is $(K,\lambda)$-stable.
\end{prop}
This is stated in \cite{D1} (Proposition 5.2.2).

  Conjecture 7.2.2 in \cite{D1} reads
\begin{conjecture}\label{conjecture_2.1.2}
If $(\Delta, K)$ is polystable, the Abreu equation $ \mc S(u)=K $
admits a solution in $\mc C_\infty(\Delta)$, where $\mc
C_\infty(\Delta)$ consists of smooth convex
 functions on $\Delta$
that are continuous on $\bar\Delta$.
\end{conjecture}
Note that the difference between $\mc C_\infty(\Delta)$ and $\mc
C^\infty(\Delta,v)$ is that the second one specifies the boundary
behavior  of the functions. On the other hand, we proved in \cite{CLS5}  the
uniform stability is a necessary condition.  Hence, related to the toric manifolds,
we state a stronger version of Conjecture \ref{conjecture_2.1.2} for Delzant
polytopes.
\begin{conjecture}\label{conjecture_2.1.3}
Let $\Delta$ be a Delzant polytope. If $(\Delta, K)$ is uniformly stable, the Abreu equation $ \mc S(u)=K $ admits a solution in
$\mc C^\infty(\Delta,v)$,
\end{conjecture}
The conjecture for cscK metric on toric surfaces  was  recently solved by Donaldson (cf.\cite{D2}).
In this paper, we solve this conjecture on toric surfaces for any edge-nonvanishing function $K.$

We need the following result proved by Donaldson.
\begin{theorem}[Donaldson\cite{D4}]\label{theorem_2.1.8}
Suppose that $\Delta$ is $(K,\lambda)$-stable. When $n=2,$ there is a constant
$\mathcal C_1>0$, depending on $\lambda$, $\Delta$ and $\|\mc
S(u)\|_{C^0}$, such that $|\max\limits_{\bar \Delta} u-\min\limits_{\bar \Delta} u|\leq \mathcal C_1$.
\end{theorem}

%%%%%%%%%%%%%%%%%%%%%%%%%%%%sect_3

\section{Some results via affine techniques}\label{sect_3}

We review the results developed in
\cite{CLS} via affine techniques.

\subsection{Calabi geometry}\label{sect_3.1}
 Let $f(x)$ be a smooth, strictly convex function defined on a
convex domain  $\Omega\subset\real^n\cong \t$.  As $f$ is strictly
convex,
$$
G:=G_f=\sum_{i,j} \frac{\partial^2 f}{\partial x_i\partial x_j}d x_i dx_j
$$
defines a Riemannian metric on $\Omega$. We call it the {\em
Calabi} metric. We recall some fundamental facts on the Riemannian
manifold $(\Omega, G)$. Let $u$ be
 the Legendre transform  of $f$ and $
 \Omega^\ast=\nabla^f(\Omega)\subset \t^\ast$.
Then it is known that $ \nabla^f: (\Omega, G_f)\to (\Omega^\ast,
G_u) $ is locally  isometric. The scalar curvature is $\mc S(f)$ or $\mc
S(u)$.

Let $ \rho=\left[\det(f_{ij})\right]^{-\frac{1}{n+2}}, $ we
introduce the following affine invariants:
\begin{equation}\label{eqn_3.1}
\Phi={\|\nabla\log\rho\|^2_G}=\frac{1}{(n+2)^2}\|\nabla\log \det(\partial^2_{ij} u)\|_{G}^2 \end{equation}
\begin{equation}\label{eqn_3.2}
4n(n-1)J =\sum f^{il}f^{jm}f^{kn}\partial^3_{ijk}f\partial^3_{lmn}f=\sum u^{il}u^{jm}u^{kn}\partial^3_{ijk}u\partial^3_{lmn}u.
\end{equation}
where $\partial^3_{ijk}f=\frac{\partial ^3 f}{\partial x_i \partial x_j \partial x_k}$ and $\partial^3_{ijk}u=\frac{\partial ^3 u}{\partial \xi_i \partial \xi_j \partial \xi_k}$. $\Phi$ is called the norm of {\em the Tchebychev vector field} and
$J$ is called {\em the Pick invariant}.
Put
\begin{equation}\label{eqn_3.3}
\Theta = J + \Phi.
\end{equation}

Consider an affine transformation
$$
\hat A: \t^\ast\times\real\to \t^\ast\times\real; \;\;\; \hat
A(\xi, \eta)=(A\xi, \lambda \eta),
$$
where $A$ is an affine transformation on $\t^\ast$. If $\lambda=1$
we call $\hat A$ the base-affine transformation. Let $\eta=u(\xi)$
be a function on $\t^\ast$. $\hat A$ induces a transformation on
$u$:
$$
u^\ast(\xi)= \lambda u(A\inv\xi).
$$
Then we have the following lemma of the {\em affine transformation
rule} for the affine invariants.
\begin{lemma}\label{lemma_3.1.1}
Let $u^\ast$ be as above, then
\begin{enumerate}
\item $\det(\partial^2_{ij} u^\ast)(\xi)=\lambda^n\det(A)^{-2}\det(\partial^2_{ij} u)(A\inv
\xi)$. \item $G_{u^\ast}(\xi)= \lambda G_u(A\inv \xi)$; \item
$\Theta_{u^\ast}(\xi)=\lambda\inv \Theta_u(A\inv \xi)$; \item $\mc
S(u^\ast)(\xi)= \lambda\inv\mc S(u)(A\inv \xi)$.
\end{enumerate}
\end{lemma}
 As a corollary,
\begin{lemma}\label{lemma_3.1.2}
$G$ and $\Theta$ are invariant with respect to the base-affine
transformation. $\Theta\cdot G$ and $\mc S\cdot G$ are invariant
with respect to affine transformations.
\end{lemma}

The following lemma was proved in \cite{CLS}.
\begin{lemma}\label{glemma_2}
  Let $u$ be a smooth, strictly convex function
defined on  $ \Omega\subset  \mathbb  R^n$ {  and $0\in \Omega$.}  Suppose that {
\begin{equation}\label{eqn_2.5}
\Theta\leq \mathsf{N}^2\quad in \quad \Omega,
\end{equation}
and the Hessian matrix $(\partial^2_{ij}u)$ satisfies $ \partial^2_{ij}u(0)=\delta_{ij}$. }
Let $\Gamma:\xi=\xi(s), s\in[0,\mff a],$ be a curve lying in $\Omega,$ starting from $\xi(0)=0$ with arc-length parameter with respect to the Calabi metric $G_{u}=\sum \partial^2_{ij}u d\xi_i d\xi_j$. Let $\lambda_{\min}$
and $\lambda_{\max}$ be the minimal and maximal eigenvalues of
$(\partial^2_{ij}u)$ along the path $\Gamma$.
Then there exists a constant  $\mff C_1$  such that
\begin{description}
\item[(i)] $ \exp\left(-\mff C_1\mff a\right)\leq
\lambda_{\min}\leq \lambda_{\max}\leq \exp\left(
\mff C_1\mff a\right), $
\item[(ii)] $\Gamma\subset D_{\mff a\exp\left(\frac{1}{2}\mff C_1\mff a\right)}(0)$, where $D_{\mff a\exp\left(\frac{1}{2}\mff C_1\mff a\right)}(0)\subset \mathbb R^{n}$ is a Euclidean disk centered at $0$ with the radius $\mff a\exp\left(\frac{1}{2}\mff C_1\mff a\right)$.
\end{description}
 \end{lemma}

\v
In  \cite{CLS}, we used the affine blow-up analysis to prove
 the following estimates. We only state the results for the Delzant
 polytopes.

\begin{theorem}\label{theorem_3.1.3} Let $u$ be a smooth strictly convex function on
a Delzant polytope  $\Delta\subset \mathbb R^2$ with $\|\mc S(u)\|_{C^{3}(\Delta)}<K_o$,
where  $\|\cdot\|_{C^3}$ denotes the Euclidean $C^3$-norm. Suppose that
for any $p \in \Delta$, $ d_u(p,\partial{\Delta})<\infty,$ and
\begin{equation}
\max_{\bar \Delta} u-\min_{\bar \Delta} u\leq \mff C_1
\end{equation}
 for some constant $\mff C_1>0.$  Then there is a constant $\mff
C_3>0$, depending only on $\Delta,\mff C_1,K_o$ such that
\begin{equation}\label{eqn_3.4}
(\Theta+|\mc S|+\mc K) d^2_u(p,\partial {\Delta})\leq \mff C_3.
\end{equation} Here $d_u(p,\partial \Delta)$ is the distance from
$p$ to $\Delta$ with respect to the Calabi metric $G_u$.
\end{theorem}

\subsection{Convergence theorems and Bernstein properties}\label{sect_3.2}
Let $\Omega \subset \mathbb{R}^2$ be a bounded convex domain. It
is well-known (see \cite{G}, p.27) that there exists a unique
ellipsoid $E$  which attains the minimum volume among all the
ellipsoids that contain $\Omega $ and that is centered at the
center of mass of $\Omega $, such that
$$2^{-\frac{3}{2}} E \subset \Omega  \subset E,$$ where
$2^{-\frac{3}{2}} E$ means the $2^{-\frac{3}{2}}$ -dilation of $E$
with respect to its center. Let $T$ be an affine transformation
such that $T(E)=D_1(0)$, the unit disk. Put $\tilde{\Omega}=
T(\Omega)$. Then
\begin{equation}\label{eqn_3.6}
2^{-\frac{3}{2}}D_1(0) \subset \tilde{\Omega} \subset
D_1(0).\end{equation} We call $T$ the normalizing transformation of
$\Omega$.
\begin{defn}\label{definition_3.2.1}
A convex domain $\Omega $ is called normalized when its center of mass is 0 and $2^{-\frac{3}{2}}D_1(0) \subset  \Omega  \subset
D_1(0).$
\end{defn}
Let $A:\real^2\to\real^2$ be an affine transformation
given by $A(\xi)=A_0(\xi)+a_0$, where $A_0$ is a linear transformation
and $a_0\in \real^2$. If there is
a constant $L>0$ such that $|a_0|\leq L$ and for any Euclidean  unit vector $v$
$$
L\inv\leq |A_0v|\leq L,
$$
we say that  $A$ is  $L$-bounded.
\begin{defn}\label{definition_3.2.2}
A convex domain $\Omega$ is called $L$-normalized if its
normalizing transformation is  $L$-bounded.
\end{defn}

The following lemma is useful to measure the normalization of a
domain.
\begin{lemma}\label{lemma_3.2.3}
 Let $\Omega \subset \real^2$ be a  convex domain. Suppose
 that there exists a pair of constants $R>r>0$ such that
 $$
D_r(0)\subset \Omega\subset D_R(0),
 $$
then $\Omega$ is $L$-normalized, where $L$ depends only on $r$
and $R$.
\end{lemma}

Let $u$ be a convex function on $\Omega$. Let
$p\in \Omega$ be a point.
Consider the set
$$
\{ \xi\in \Omega| u(\xi)\leq u(p)+\nabla
u(p)\cdot(\xi-p)+\sigma\}.
$$
If it is compact in $\Omega$, we call it {\em a section of $u$ at
$p$ with height $\sigma$}  and denote it by $S_u(p,\sigma)$.

\v\n
Denote by $\mathcal{F}(\Omega,C)$ the class of convex functions
defined on $\Omega$ such that
$$ \inf_{\Omega} {u } = 0,\;\;\;
u= C\;\;on\;\;\partial \Omega,$$ and
$$
\mc F(\Omega,C;K_o)=\{u\in \mc F(\Omega,C)| |\mc S(u)|\leq K_o\}.
$$
We will assume that  $\Omega$ is normalized in this subsection.

The main result of this subsection is the following.
\begin{prop}\label{proposition_3.2.4}
Let $\Omega\subset \mathbb R^2$ be a  normalized domain. Let $u\in \mc F(\Omega,1;
K_o)$ and $p^o$ be its minimal point, that is, $u(p^{o})=0$. Then
\begin{enumerate}
\item[(i)] there are two positive  constants $\mff s$ and $\mff C_2$ such that $d_E(p^o,\partial\Omega)>\mff s$ and in
$D_{\mff s}(p^o)$
$$
\|u\|_{C^{3,\alpha}}\leq \mff C_2
$$
for any $\alpha\in (0,1)$;
\item[(ii)] there is a constant
$ \delta\in (0,1),$ such that $S_u(p^o,\delta)\subset D_{\mff s}
(p^o). $ \item[(iii)] there exists a constant $b>0$ such that
$S_u(p^o,\delta)\subset B_b(p^o)$.
\end{enumerate}
\end{prop}
In the statement,
all the constants  only depend on $K_o$; $D, B$ are
   disks
with respect to the Euclidean metric and the Calabi metric $G_u$ respectively;
$d_E$ is the Euclidean distance function. Here and later we denote  $\|\cdot\|_{C^{3,\alpha}}$ the Euclidean $C^{3,\alpha}$-norm.

Furthermore, if $u$ is smooth, then for any $k\in \mathbb Z_{\geq 0}$
$$
\|u\|_{C^{k+3,\alpha}(D_{\mff s}(p^o))}\leq \mff C_2'
$$
where $\mff C_2'$ depends on the $C^{k}$-norm of $\mc S(u)$.

\v This can be restated as a convergence theorem.
\begin{theorem}\label{theorem_3.2.5}
Let $\Omega\subset \mathbb R^2$ be a  normalized domain. Let $u\indexm\in \mc
F(\Omega,1; K_o)$ be a sequence of functions and $p^o_k$ be
the minimal point of $u\indexm$. Then there exists a subsequence
of functions, still denoted by
$u\indexm$, locally uniformly converging to a function $u_\infty$ in $\Omega$, and
$p^o_k$ converging to $p^o_\infty$; satisfying:
\begin{enumerate}
\item[(i)] there are two positive constants $\mff s$ and $\mff C_2$ independent of k such that $d_E(p_k^o,\partial\Omega)>\mff s$ and in
$D_{\mff s}(p^o_\infty)$
$$
\|u\indexm\|_{C^{3,\alpha}}\leq \mff C_2
$$
for any $\alpha\in (0,1)$;
 \item[(ii)] there is a
constant $ \delta\in (0,1),$  independent of k, such that $S_{u\indexm}(p^o_k,
\delta)\subset D_{\mff s} (p^o_\infty). $ \item[(iii)] there
exists a constant $b>0$  independent of k  such that
$S_{u\indexm}(p^o_k,\delta)\subset B_b(p_k^o)$.
\end{enumerate}
\end{theorem}
(i) implies that in
$D_{\mff s}(p^o_\infty)$, $u\indexm$ $C^{3}$-converges
to $u_\infty$.
Furthermore, if $u\indexm$ is smooth and the  $C^{k}$-norms
of $\mc S(u\indexm)$ are uniformly bounded, then $u\indexm$
$C^{k+3,\alpha}$ converges to $u_\infty$ in
$D_{\mff s}(p^o_\infty)$.

\subsection{Interior regularities and estimate of $\mc K$ near divisors}
\label{sect_3.4}
\def\fkz{\mathfrak{z}}
Let $\Delta\subset \mathbb R^2$ be a Delzant  polytope.
In \cite{CLS}, we prove the following regularity theorem.
\begin{theorem}\label{theorem_3.4.1}
Let $\CHART$ be a chart of either $\CHART_\Delta, \CHART_\ell$ or $\CHART_\vartheta$.
 Let $z_o\in\CHART$ and
 $B_a(z_o)$ be a geodesic ball in $\CHART$.
 Suppose that there is a
 constant $C_1$ such that $f(z_o)=0,\;\;|\nabla f|(z_o)\leq C_1,\;\;$ and
$$\mc K(f) \leq C_1,\;\;\;\;   W \leq   C_1,\;\;\; |z|\leq C_1 $$
in  $B_a(z_o)$.
 Then  there is a constant $a_1>0$, depending on $a$ and $C_1$,
such that  $D_{{2a_1}}(z_o)\subset B_{\frac{a}{2}}(z_o),$ and for any $k\geq 0,$ 
such that
$$\|f\|_{ C^{k+3,\alpha}(D_{  a_1}(z_o ))}\leq C(a,C_1, \|\mc S(u)\|_{C^{k}}).$$
\end{theorem}

\v One of the main results in \cite{CLS} that is developed from the
affine technique is the following. Here $\mathcal K, W$ and $\Psi$ are introduced in
\S\ref{sect_1.2}.
 \begin{theorem}\label{theorem_3.4.2}
Let $u\in \mc C^\infty(\Delta,v)$. Let $\fkz_\ast$ be a
point on a divisor $Z_\ell$
for some $\ell$. Choose a coordinate system $(\xi_1,\xi_2)$ such that $\ell=\{\xi|\xi_1=0\}.$ Let $p\in \ell$ and $D_b(p)\cap \bar \Delta$ be an Euclidean half-ball such that its intersects with $\partial\Delta$ lies in the
interior of $\ell$. Let $B_a(\fkz_\ast)$ be a geodesic
ball satisfying  $\tau_f(B_a(\fkz_\ast))\subset D_b(p)$.
 Suppose that
\begin{eqnarray}
 |\mc S(u)|\geq \delta>0, \;in\; \;D_b(p)\cap \bar \Delta,\\
 \|\mc S(u)\|_{C^3(\bar \Delta)}\leq\mff N ,  \nonumber\\
 \partial^2_{22}h|_{\ell\cap D_b(p) }\geq \mff N \inv \nonumber
\end{eqnarray}
    where $h=u|_{\ell}$ and $\|.\|_{C^3(\Delta)}$ denotes the Euclidean  $C^{3}$-norm.
Then there is a
 constant $\mff C_3>0$,
 depending only on $a,\delta,\mff N$ and  $D_b(p)$, such that
\begin{equation}\label{est_eq_1.2}
 \frac{W^{\frac{1}{2}}}{\max\limits_
 {B_a(\fkz_\ast )} W^{\frac{1}{2}}}
\left( \mc K+\|\nabla \log |\mc S|\|_f^2 +\Psi\right)(\fkz)a^2 \leq \mff C_3,\;\;\;\;\;\;\;\forall \fkz\in B_{a/2}(\fkz_\ast)
\end{equation}
where $W=\detf$.
\end{theorem}

%%%%%%%%%%%%%%%%%%%%%%%%%%%%%%%%sect_4

\section{Estimates of the determinant}
\label{sect_4}

%As stated in Theorem \ref{theorem_2.2.3}, the interior regularity in $\Delta$ has been established.
 In this section, we will explore the dependence of
 some estimates of $\det(\partial^2_{ij}u)$ on $d_E(\cdot,\partial\Delta)$.
The results in \S\ref{sect_4.1}   hold for any $n$.
\subsection{The lower bound of the determinant}
\label{sect_4.1}
The following  lemma can be found in \cite{D3}.
\begin{lemma}\label{lemma_4.1} Suppose that $u\in \mc C^\infty(\Delta,v;K_o)$. Then
\begin{description}
\item[(1)] $\det (\partial^2_{ij}u)\geq \mff C_4$ everywhere in $\Delta$,
where $\mff C_4 = ({4n\inv K_o
 \mathrm{diam}(\Delta)^2})^{-n}.$
\item[(2)] For any $\delta\in (0,1)$ there is a constant $\mff
C_\delta>0$, depending only on $n$ and $\delta$,  such that
$$\det(\partial^2_{ij}u)(p)\geq\mff  C_\delta d_E(p,\partial
\Delta)^{-\delta}.$$
\end{description}
\end{lemma}
Here  we denote $\partial^2_{ij}u=\frac{\partial^2u}{\partial \xi_i\partial \xi_j}.$
In the following  we derive a stronger estimate for $\det(\partial^2_{ij}u)(p)$.

\begin{lemma} \label{lemma_4.2}
Let $u\in \mc C^\infty(\Delta,v;K_o)$.   Then there is a constant $\mff
C_5>0$, depending only on $\Delta$ and $K_o$ such that for any
$\xi\in \Delta$
$$
\det(\partial^2_{ij}u)(\xi)\geq \frac{\mff C_5}{d_E(\xi,\partial \Delta)}.
$$
\end{lemma}
 {\bf Proof.} Let $p\in \Delta$ be a point and $F$ be a facet such that $d_E(p,\partial \Delta)=d_E(p,q),$ $q\in F$.
We choose a new coordinate system on $\t^\ast$ such that: {\em (i)
 $F$ is on the $\{\xi_1=0\}$-plane; (ii) $\xi(q)=0$ ;  (iii)
$\xi_1(\Delta)\geq 0$. }

By (2) of Lemma \ref{lemma_4.1}, we already have
\begin{equation}\label{eqn_4.1} \det(\partial^2_{ij}u)\geq
C_0\xi_1^{-(1-\frac{1}{n})}.\end{equation} Consider the
function
$$v'= \xi_1^{\alpha}\left(C + \sum_{j=2}^{n}
\xi_j^2\right) - a \xi_1,$$ where $a>0$, $\alpha >1$ and $C>0$ are
constants to be determined. We may choose $a$ large such that
$v'\leq 0$ on $\Delta$. For any point $\xi$ we may assume that
$\xi = (\xi_1,\xi_2,0,...,0)$. By a direct calculation we have
$$v'_{11}= \alpha (\alpha-1)\xi_1^{\alpha-2}(C+
\xi_2^2),\;\;\; v'_{12}= 2\alpha \xi_2\xi_1^{\alpha-1}, \;\;\;
v'_{ii}= 2\xi_1^{\alpha}\;\;\;i\geq 2,$$
$$\det(\partial^2_{ij}v')=
2^{n-1}\left[\alpha (\alpha-1)(C+ \xi_2^2) -
2\alpha^2\xi_2^2\right]\xi_1^{n\alpha-2}.$$
Set $\alpha=1+\frac{1}{n^2}$. Then for large $C$, it is easy to see that $v'$ is strictly convex in $\Delta $
and
\begin{equation}\label{eqn_4.2}
\det(\partial^2_{ij}v')\geq C_1\xi_1^{n\alpha-2}.
\end{equation}

 Consider
the function
$$F = w + C_5v',$$
where $w\inv=\det(\partial^2_{ij}u)$. As $w$ vanishes  on the boundary of
$\Delta$, we have $F\leq 0$ on $\partial \Delta$. Then
\begin{eqnarray*}
\sum U^{ij}\partial^2_{ij}F& =& -K +C_5\det(\partial^2_{ij}u)\sum u^{ij}\partial^2_{ij}v'\\
& \geq&
-K + nC_5 \det(\partial^2_{ij}u)^{1-1/n}\det(\partial^2_{ij}v')^{1/n}\\
&\geq& -K + nC_5
C_{0}^{1-\frac{1}{n}}\xi_1^{-(1-\frac{1}{n})^2}C_1^{\frac{1}{n}}\xi_1^{\alpha-\frac{2}{n}}
\\
&=&-K+nC_5C_{0}^{1-\frac{1}{n}}C_1^{\frac{1}{n}}.
\end{eqnarray*}
where $C_{0}$ and $C_1$ are constants in \eqref{eqn_4.1} and \eqref{eqn_4.2}, $(u^{ij})$ denotes the inverse matrix of the matrix $(\partial^2_{ij}u)$. Choose $C_5$ such that
$\sum U^{ij}\partial^2_{ij}F> 0$. So by the maximum principle we have $w \leq C_5|v'|\leq aC_5\xi_1$. It
follows that % for any facet $F$,
 $\det(\partial^2_{ij}u)(\xi)\geq {\mff C_5}{\xi_1}\inv $ for some
constant $\mff C_5>0$ independent of $p.$
  $\blacksquare$

\subsection{ The upper bound of the determinant }\label{sect_4.2}
 We need the
following two lemmas.
\begin{lemma}\label{lemma_4.2.3}
Let $\Delta  \subset \mathbb R^n$ be a Delzant ploytope and $p\in \Delta$ be a point with
$ \xi(p)=0.$ Let $u_k\in \mc C^\infty(\Delta)$ be a sequence of convex functions such that
$$ u_k(0)=0,\;\;\nabla
u_k(0)=0,$$ and  $u_k$ locally uniformly $C^{2 }$-converges to a strictly convex function $u_\infty$ defined in $\Delta.$    Then there are two constants $d,C_1>0$ independent of
$k$ such that
$$\frac{\sum\left(\frac{\partial u_k}{\partial \xi_i}\right)^2}{(d+ f_k)^2}\leq C_1$$
where   $f_k$ is the Legendre function of $u_k$ (cf. Section \ref{sect_1.3}).
\end{lemma}
{\bf Proof.}  Obviously $f_k(0)=0,\;\;\nabla f_k(0)=0,\;\;f_k\geq 0$
and $f_k$ uniformly
$C^{2 }$-converges to a strictly convex function $f_\infty$ in
$D_\epsilon(0)$ for some $\epsilon>0$, in particular,
$$f_k|_{\partial D_\epsilon(0)}\geq \delta $$
for some $\delta>0$ independent of $k$.
 Let $h(x)=\frac{\delta|x|}{\epsilon}.$ Using the convexity of $f_k$
 one can check that in $\mathbb R^n \setminus D_\epsilon (0),$
 $$ f_k(x)\geq h(x) .$$
 Then in $\mathbb R^n$
 $$(f_k+\delta)^2\frac{\epsilon^2}{\delta^2} \geq
 \sum_{i} x_i^2 =  \sum_{i}\left(\frac{\partial u_k}{\partial \xi_i}\right)^2  .\;\;\;\;  q.e.d.$$

 The following lemma is proved in \cite{CLS1} (cf. Corollary 2.6):

\begin{lemma}\label{lemma_4.2.4} Let $\Delta\subset \mathbb R^2$ be a Delzant ploytope.
 Suppose that $u\in \mc
C^\infty(\Delta,v;K_o)$, and suppose that there are two constants
$b,d>0$ such that
\begin{equation}\label{eqn_2.0.a}
\frac{\sum \left(\frac{\partial u}{\partial \xi_k}\right)^2}{(d+
f)^2}\leq b,\;\;\;\; d+f\geq 1
\end{equation} where  $f$ is the Legendre function of $u.$
  Then there is a constant $b_0>0$ depending only on $K_o$ and $\Delta$ such that
$$\frac{\det(\partial^2_{ij}u)}{(d+f)^4}(p) \leq \frac{b_0}{ d_E(p, \partial
\Delta)^{4}}.$$
\end{lemma}

Using Lemma \ref{lemma_4.2.3} and Lemma \ref{lemma_4.2.4} we can obtain the upper bound estimates for $\det(\partial^2_{ij}u).$

\begin{lemma}\label{lemma_4.2.5} Let $\Delta\subset \mathbb R^2$ be a Delzant ploytope. Suppose that $u_k\in \mc
C^\infty(\Delta,v;K_o)$ and  $u_k$ locally uniformly $C^{2}$-converges
to a strictly convex function $u_\infty$ in $\Delta.$ And suppose that
$$\max_{\bar \Delta} |u_k|\leq \mc C_1,$$
 for some constant $\mc C_1>0$ independent of $k$. Denote    $d_E(p,\partial \Delta)$ by   the Euclidean distance from $p$ to the boundary $\partial\Delta $. Then
 there is a constant  $\mff C_6 >0,$  independent of $k$,  such that for any $p\in \Delta$
  $$\log\det(\partial^2_{ij} u_k)(p)\leq \mff C_6- \mff  C_6  \log d_E(p,\partial \Delta).$$
\end{lemma}
{\bf Proof.}
Since $u$ is convex, we have for any $p\in \Delta$
\begin{equation}\label{eqn_a5.3}
\left|\frac{\partial u}{\partial \xi_i}(p)\right|\leq \frac{2\mc C_1}{d_E(p,\partial \Delta)}.
\end{equation}
  Again by the convexity and \eqref{eqn_a5.3} we obtain that
for any point $p\in \Delta$
\begin{equation}\label{eqn_4.2.5}f(\nabla_u(p))=\sum \frac{\partial u}{\partial \xi_i}\xi_i -u \leq \frac{4C_2\mc C_1}{d_E(p,\partial \Delta)},\end{equation}
 where $C_2>0$ is a constant depending only on $\Delta.$
  From Lemma \ref{lemma_4.2.3}, Lemma \ref{lemma_4.2.4} and  \eqref{eqn_4.2.5}  we
conclude that
\begin{equation}
\det(\partial^2_{ij}u)\leq C_3(d_E(p,\partial \Delta))^{-8},\;\;\;\;
\end{equation}
where $C_3>0$ is a constant.   \;\;\;\;  q.e.d.

%%%%%%%%%%%%%%%%%%%%%%%%%%%%%%%%%%sect_5

\section{Estimates of Riemannian distances on $\partial\Delta$}
\label{sect_5}
Let $\Delta\subset \mathbb R^2$ be a Delzant ploytope. Let $\ell$ be an edge of $\Delta$ and
$\ell^\circ$ be the interior of $\ell$.  Let $\xi^{(\ell)}\in \ell^{\circ}$. For simplicity, we fix a coordinate system on
$\t^\ast$ such that { (i) $\ell$ is on the $\xi_2$-axis; (ii) $\xi^{(\ell)}=0$;
(iii) $\Delta\subset \halfplane$.}

Define $\ell_{c,d}=\{(0,\xi_2)|c\leq \xi_2\leq d\}\subset \ell^o$.

Let $u\indexm \in \mc C^\infty(\Delta,v;K_o)$ be a sequence of
functions with $\mc S(u\indexm)=K\indexm$. Suppose that
\begin{enumerate}
\item \begin{equation}\label{eqnc_5.1}
\left|\max_{\bar\Delta} u_k-\min_{\bar\Delta} u_k\right| \leq \mc C_1
\end{equation}  for some constant $\mc C_1$ independent of $k$,
\item $K\indexm$   $C^3$-converges to $K$ on $\bar \Delta$, and
\item  $u\indexm$  locally   $C^6$-converges
in $\Delta$ to a strictly convex function $u_\infty$. $u_\infty$ can be naturally continuously  extended to be defined on $\bar \Delta.$
\end{enumerate}

\subsection{$C^0$-convergence}\label{sect_5.1}
Denote by $h\indexm$
the restriction of $u\indexm$ to $\ell$. Then $h\indexm$ locally
uniformly converges to a convex function $h $ in $\ell$.
Obviously, $u_\infty|_{\ell^\circ} \leq h $.
In this
 subsection we prove that "$\leq$" is indeed "$=$".
In fact, we have the following Proposition.
\begin{prop}\label{proposition_5.1.1}
For $q\in \ell^\circ$,
$u_\infty(q) = h(q)$.
\end{prop}
{\bf Proof.} For simplicity we assume that
\begin{equation}\label{eqn_5.1a}
\xi(q)=0.\end{equation}    If this proposition  is not true, then
$u_\infty(0) < h(0)$. Without loss of generality  we can assume that $
\ell_{-\half,\half}=\{\xi|\xi_1=0,|\xi_2|\leq \half\} \subset
\ell^o  $ and  for any point
$p\in \ell_{-\half,\half}$
$$u_\infty(p) +\half < h(p),\;\;\; u_\infty(0)=0 .$$
By assumption we have
\begin{equation}\label{eqna_5.1}
\lim_{k\to\infty}\|K_k-K\|_{C^3(\bar \Delta)}=0.
\end{equation}
 For any $K_k,$ consider the functional
$$ \mc F_{K_k}(u)= -\int_{\Delta}\log \det(\partial^2_{ij}u)d\mu + \mc L_{K_k}(u),$$
defined in $\mathcal C_\infty(\Delta),$
where $\mc L_{K_k}$ is the linear functional
$$\mc L_{K_k}(u) = \int_{\partial \Delta}u d\sigma - \int_{\Delta}K_ku
d\mu.$$ Here $d\sigma$ and $d\mu$ are as in the section \ref{sect_2.1}. Since $\mc S(u\indexm)=K\indexm$, by a result of Donaldson  $u_k$ is an absolute minimizer for
$\mc F_{K_k}$ in  $\mathcal{C}_{\infty}(\Delta)$ (cf. \cite{D1}).
By \eqref{eqnc_5.1} and $u_\infty(0)=0$
\begin{equation}\label{eqna_5.2}  |u_k|_{L^\infty(\bar \Delta)}<\mc C_1  \end{equation}
as $k$ large enough.
For any positive constant $\delta<1$, denote
 $$\Delta_\delta =\{p\in \Delta\;|\; d_E(p,\partial \Delta)\geq \delta\}.$$
 By Lemma \ref{lemma_4.1} and Lemma \ref{lemma_4.2.5} we have for any $k$
  $$\left|\int_{\Delta\setminus \Delta_\delta}\log \det(\partial^2_{ij}u_k)d\mu\right|\leq C_2\sqrt{\delta}$$
  as $\delta$ is small enough,
  where $C_2$ is a constant independent of $k$ and $\delta$. Combining this and that  $u_k$ locally uniformly $C^3$-converges to $u_\infty$ in $\Delta$ we have
$$\lim_{k\to\infty} \int_{\Delta}\log \det(\partial^2_{ij}u_k)d\mu=
\int_{\Delta}\log \det(\partial^2_{ij}u_\infty)d\mu.\;$$
By \eqref{eqna_5.1}
 and \eqref{eqna_5.2} we have
$$
\lim_{k\to\infty}\int_{\Delta}K_ku_kd\mu= \int_{\Delta}K u_\infty d\mu
$$
and \begin{equation*}\lim_{k\to\infty}\int_{\partial \Delta}u_k
d\sigma- \int_{\partial \Delta}u_\infty d\sigma=\int_{\partial
\Delta}h d\sigma- \int_{\partial \Delta}u_\infty d\sigma \geq
\int_{\ell_{-\half,\half}}(h-u_\infty) d\sigma\geq \half.
\end{equation*}
We conclude that
\begin{equation}\label{eqn_5.2}
\mc F_{K}(u_\infty)\leq\lim_{k\to\infty}\mc F_{K_k}(u_k) -\frac{1}{2}.
\end{equation}
Hence
\begin{equation}
\mc F_{K_k}(u_\infty)=\mc  F_{K}(u_\infty)-\int_{\Delta }( K_k-K) u_\infty d\mu \leq \mc F_{K_k}(u_k)-\frac{1}{4},\end{equation}
 as $k$ large enough, where we used \eqref{eqna_5.1}
 and \eqref{eqna_5.2} in the last inequality.
This contradicts  $\mc F_{K_k}(u_k)=\inf\limits_{u\in \mathcal C_\infty(\Delta)}\mc F_{K_{k}}(u)$.
$\blacksquare$

\begin{corollary}\label{corollary_uniform_control}
 Suppose that $\ell_{c-2\epsilon_{o},d+2\epsilon_{o}}\subset \ell $ for some $\epsilon_{o}>0.$ For any $\epsilon>0,$ there is a constant $\delta>0$, such that for any $ q\in \ell_{c,d}$ and $p\in E_{c,d}^{\delta}:=\{(\xi_1,\xi_2)\in \Delta | \xi_{1}\leq \delta, c\leq \xi_{2}\leq d\}$ with $d_E(p,q)\leq \delta,$
\begin{equation}\label{c_o_uniform_controll}
|u_{k}(p)-u_{k}(q)|\leq \epsilon,
\end{equation}
when $k$ is big enough.
\end{corollary}
{\bf Proof.} Let $\delta_{0}$ be a constant such that $\overline{E_{c- \epsilon_{o},d+ \epsilon_{o}}^{2\delta_{0}}}\cap \partial \Delta \subset \ell.$
By the convexity of $u_{k}$, for any $p\in E_{c,d}^{\delta_{0}}$, we have
\begin{equation}\label{jr_eqn5.7}
\partial_{1}u_{k}(p)< \mathcal C_{1}\delta_{0}\inv,\;\;\;\;\;\;|\partial_{2}u_{k}(p)|\leq   \mathcal C_{1} \epsilon_{o}\inv.
 \end{equation}
Without loss of generality  we can assume that $\max_{E_{c,d}^{\delta_{0}}}\partial_1 u_k\leq 0$ (otherwise, we can use $u-\mathcal C_{1}\delta_{0}\inv \xi_1$ instead of $u,$   and use the same argument).
 Let $0<\delta_{1}\leq  \delta_{0} $ be a constant such that
\begin{equation}\label{jr_eqn5.8}
\max_{c\leq b\leq d} |  u_{\infty}(0,b)-  u_{\infty}(\delta_{1},b)|\leq  \epsilon/8.
\end{equation}
By \eqref{jr_eqn5.7}  we conclude that $u_{k}$ is uniform continuous in the $\xi_{2}$-direction. As $\max\limits_{E_{c,d}^{\delta_{0}}}\partial_{1}   u_{k}<0$ we have $u_{k}(\delta_1,b)<  u_{k}(a,b)<  u_{k}(0,b)$ for any $0<a<\delta_1.$  It   suffices to prove that for any $c\leq  b\leq d$
   $$|u_{k}(0,b)-u_{k}(\delta_1,b)|\leq \epsilon.$$
By the convergence of $h_{k}$ and the convergence of $  u_{k}$   we have
\begin{equation}\label{jr_eqn5.9}
\max_{c\leq b\leq d}| u_{k}(\delta_{1},b)- u_{\infty}(\delta_{1},b)|\leq \frac{\epsilon}{8},\;\;\;\;\;\;
\max_{c\leq b\leq d}| h_{k}(0,b)-h(0,b)|\leq \frac{\epsilon}{8},
\end{equation}
when $k$ large enough.   By Proposition \ref{proposition_5.1.1} we have
\begin{align*}& |   u_{k}(0,b)-  u_{k}(\delta_1,b)|  \\
&\leq  |   u_{\infty}(0,b)-  u_{\infty}(\delta_1,b)|+ |   u_{\infty}(\delta_1,b)-  u_{k}(\delta_1,b)|+ |   h_{\infty}(0,b)-  h_{k}(0,b)| .
\end{align*}
Then the Corollary follows from  \eqref{jr_eqn5.8} and \eqref{jr_eqn5.9}.  $\Box$

\subsection{Monge-Amp\`ere measure on the boundary}
\label{sect_5.2}

\begin{lemma}\label{lemma_5.2.1}
Let $u\in \mc C^\infty(\Delta,v;K_o)$ and $h=u|_{\ell}$. There is
a constant
 $\mff C_7>0$, depending on $d_E(\ell_{c,d},\partial\ell)$, such that on
$\ell_{c,d}$, $ \partial^2_{22}h\geq \mff C_7.$\end{lemma}
 {\bf Proof.} By the boundary behavior of $u$,
we know that $\frac{\partial u^{11}}{\partial \xi_1} = 1$  on $\ell$ (cf. \cite{D4}). Consider a
small neighborhood of $\ell_{c,d}$ which depends on $u,$ such that $\frac{1}{2}\leq
\frac{\partial u^{11}}{\partial \xi_1}\leq 2$.  By integrating we have, in this neighborhood,
$$\frac{1}{2}\xi_1\leq u^{11}\leq 2\xi_1.$$
Then, by Lemma \ref{lemma_4.2},
$\partial^2_{22}u=\det(\partial^2_{ij}u)\cdot u^{11}\geq \mff C_7.$
$\blacksquare$

\v
For any $\delta>0$ denote $$
L^{\delta}_{c}=\{(\xi_1, c)|0\leq \xi_1\leq \delta\},\;\;\;  L^{\delta}_{d}=\{(\xi_1,
d)|0\leq \xi_1\leq \delta\}.
$$

\begin{lemma} \label{lemma_grad_control}
Suppose that $\ell_{c-\epsilon_{o},d+\epsilon_{o}}\subset \ell $ for some $\epsilon_{o}>0.$  Then there exists $\delta>0$ independent of $k,$ such that
\begin{equation}\label{add_grad_control}
\max_{L^{\delta}_{c}}\partial_{2}u_{k}< \min_{L^{\delta}_{d}}\partial_{2}u_{k},
\end{equation}
when $k$ is big enough.
\end{lemma}
{\bf Proof.} Since \eqref{add_grad_control} is invariant under adding a linear  function, we can assume that
$$ \partial_{2} u_{k}(0,e)=0,\;\;\;\;\; u_{k}(0,e)=0, $$
where $e=(c+d)/2.$ Using Lemma \ref{lemma_5.2.1}, a direct integration gives us
$$ u_{k}(0,c)\geq\delta_{0},\;\;\;\; u_{k}(0,d)\geq\delta_{0},\;\;\;\;  $$
for some $\delta_{0}>0$ depends only on $\mff C_7$ and $d-c.$
% It follows that $$h(0,c)\geq\delta_{0},\;\;\;\; h(0,d)\geq\delta_{0}.$$
 Applying Corollary \ref{corollary_uniform_control} with $\epsilon=\frac{\delta_{0}}{10},$ there is a constant $\delta>0$, such that for any $0 \leq a\leq \delta$ and $c\leq b \leq d$
\begin{equation}
|u_{k}(a,b)-u_{k}(0,b)|\leq \frac{\delta_{0}}{10}
\end{equation}
when $k$ is big enough. Then
$$u_{k}|_{L_{c}^{\delta}}\geq \frac{9\delta_{0}}{10},\;\;\;\; u_{k}|_{L_{d}^{\delta}}\geq \frac{9\delta_{0}}{10},\;\;\;\; -\frac{\delta_{0}}{10}\leq u_{k}|_{L^{\delta}_{e}}\leq \frac{\delta_{0}}{10}.$$
By the convexity of $u_{k},$ we have
 $\partial _{2} u_{k}|_{L_{c}^{\delta}}<0$ and $\partial _{2} u_{k}|_{L_{d}^{\delta}}>0.$  $\;\;\;\;\Box$

\subsection{Some lemmas}
Let $u\in \mc C^\infty(\halfplane,v_{\halfplane}; K_o)$. Let $p^\circ$
be a point such that $d(p^\circ, \t_2^\ast)=1$, where $\t^\ast_2=\partial\halfplane$.
 By adding a linear function we normalize $u$ such that $p^\circ$ is the minimal point of $u$; i.e.,
 \begin{equation}\label{aeqnc_7.4}
  u(p^\circ)=\inf_{\halfplane} u.
 \end{equation}
Let $\check p$ be the minimal point of $u$ on $\t_2^\ast$, the boundary of $\halfplane$.
By adding some constant to $u$, we require that
\begin{equation}\label{aeqnc_7.5}
u(\check p)=0.
\end{equation} By a coordinate translation we can assume that \begin{equation}\label{aeqnc_7.6} \xi(\check p)=0.\end{equation}
 We call $(u,p^\circ,\check p)$   a  normalized triple,
if $u$ satisfies \eqref{aeqnc_7.4}, \eqref{aeqnc_7.5}, \eqref{aeqnc_7.6}  and
$$d(p^\circ, \t_2^\ast)=1.$$

\begin{lemma}
Let $(u\indexm, p_k^\circ,\check p\indexm)$ be a sequence  of normalized triples  with
\begin{equation}\label{eqn_a}
\lim_{k\to \infty} \max|\mc S(u\indexm)|
=0,\;\;\;\;\Theta_{u\indexm} d^2_{u\indexm}(p,\t_2^\ast)\leq \mff C_3,
\end{equation}
then there is a constant $C_o>0$ such that
\begin{equation}\label{geqn_9.1}
C_{o}\inv \leq |u\indexm(p^{\circ}_{k})|\leq C_o.
\end{equation}
\end{lemma}
The proof of this Lemma is the same as Lemma 7.6 in \cite{CLS}.

\v
Based on this, we prove
the following lemma in this subsection.
\begin{lemma}\label{glemma_3}
Let $(u\indexm, p_k^\circ,\check p\indexm)$
 be a sequence of normalized triple  with \eqref{eqn_a} and $\partial^2_{ij}u\indexm(p_k^\circ)=\delta_{ij}$.
  Then by choosing a subsequence we have
\begin{enumerate}
\item[(1)]
  $p_k^\circ$ converges to a point $p^\circ_\infty,$ %and $u_k$ locally uniformly $C^3$-converges to a strictly convex function $u_\infty,$ in particular,
  there exists a constant $a>0$ such that $u_k$ uniformly $C^3$-converges to
$u_\infty $  in a Euclidean ball $D_{a}(p^{\circ}_\infty)$;
\item[(2)] there exist two constants $0<\tau<1$ and $C_1>0$ independent of k such that
   \begin{eqnarray}
    && \label{geqn_2.5} \max_{S_{h }(\check p,1)}|\xi_2|\leq \frac{C_1}{2},\\
     && \label{geqn_2.6}\max_{S_{h }(\check p,1)}\partial_2 u \geq C_1\inv,\;\;\;\min_{S_{h }(\check p,1)}\partial_2 u \leq -C_1\inv,\;\;\;\\
  && \label{geqn_2.7}|\nabla u |\leq
     C_1\inv,\;\;\;in\; B_{\tau}(p^\circ )
    \end{eqnarray}
    as k large enough, where $h =u |_{\t_2^\ast}$ and $S_{h }(\check p,1)=\{\xi\in \t_2^\ast\;|\; h \leq 1\}.$
\end{enumerate}
\end{lemma}
{\bf Proof.} Let $u$ be a function of $u_k$. By a coordinate translation
$\xi^{\star}=\xi-\xi(p^\circ) $ we have  $\xi^{\star}( {p}^\circ)=0.$

 By \eqref{eqn_a} and $d(p^\circ, \t_2 ^\ast)=1,$ we
have
$$\Theta\leq 16\mff C_3,\;\;in\;\; B_{\frac{3}{4}}(  p ^\circ) .$$ Using
Lemma \ref{glemma_2} we obtain
$$C_{2}\inv \leq \partial^2_{ij}u\leq C_{2},\;\;\;\;|\partial^3_{ijk}u|\leq C_{2} .$$ It follows from \eqref{aeqnc_7.4} and \eqref{geqn_9.1} that $\|u\|_{C^{3}(B_{\frac{3}{4}}(   p ^\circ))} \leq C_{1}.$ Then $U^{ij}\in  C^{1}(B_{\frac{3}{4}}(   p ^\circ)) .$
 Following from the standard elliptic regularity theory of the equations
 \[ \sum U^{ij}\partial^2_{ij}w=-K,\;\;\;\;\; U^{ij}\partial^2_{ij}(\partial_{k}u)=\partial_{k}w^{-1} \]
  we have
$\|u\|_{W^{4,p}(B_{\frac{3}{4}}(   p ^\circ))}\leq C $. By the Sobolev embedding theorem
\[\|u\|_{C^{3,\alpha}(B_{\frac{1}{2}}(   p ^\circ))}\leq C_2  \|u\|_{W^{4,p}(B_{\frac{3}{4}}(   p ^\circ))}.
\]
for some positive constant $C_2$ independent of $k$.
Then by Ascoli Theorem and choosing a subsequence  we conclude that $u\indexm$ uniformly $C^3$-converges to a strictly convex function $u_\infty$ in $D_a(0)$ for some constant $a>0.$ In particular, there is a positive  constant    $\epsilon$ such that
\begin{equation}\label{geqn_2.8}
S_{u\indexm}(0,\epsilon)\subset D_a(0),\;\;\; \left|\nabla_{e_{r}} u\indexm\right|(p)\geq \frac{\epsilon}{2a},\;\;\;\;\forall p\in \partial D_a(0).
\end{equation}
where $e_r$ is a unit vector parallel to $\overline {p^{\circ}p}.$ Consider the function $$\Lambda(\xi^\star) =\frac{\epsilon}{2a}| \xi^\star|-C,$$
 where $C$ is the constant in  \eqref{geqn_9.1}.
Since $\Lambda(0)<u(0)$ and $\Lambda(p)\leq u(p) $ for any
$p\in \partial D_a(0),$ by  \eqref{geqn_2.8} and the convexity of $u$ we have
\begin{equation}\label{geqn_2.11}
\Lambda(q)\leq u(q),\;\;\;\;\;\; \forall \; q\in \mathsf{ h}^\ast
\setminus D_a(0).
\end{equation}

Then by $ \t^\ast_{2 }\subset  \mathsf{  h}^\ast
\setminus D_a(0),$ we have
 $
S_h( {\check p},1)\subset S_{\Lambda}({  p}^\circ,C+1)\cap  \t_{2}^\ast.
$
In particular
 \begin{equation}
 \label{geqn_c.8}
\max_{S_h( {\check p},1)} |\xi^\star| \leq  \max_{S_{\Lambda}({p}^\circ,C+1)} | \xi^{\star}| \leq \frac{2a}{\epsilon}(C_1+1). \end{equation}
Combing this, $\xi^{\star}(\check{p})=-\xi(p^{\circ})$ and $ \check{p}\in S_h( {\check p},1)$, we prove \eqref{geqn_2.5} and
\begin{equation}\label{alpha_beta_c.8}
|\xi(p^{\circ})|\leq \frac{2a}{\epsilon}(C_1+1).
\end{equation}
 (1) follows from \eqref{alpha_beta_c.8} and  the convergence of $u\indexm$ in the coordinates of $\xi^\star$.
\eqref{geqn_2.7} follows from (1) and the convexity of $u_\infty$.
Then by the convexity of $u$ and \eqref{geqn_2.5} we obtain \eqref{geqn_2.6}.
   $\Box$

\subsection{Lower bounds of Riemannian distances inside edges}
\label{sect_5.3}
 Let $p=(0,c)$ and
$q=(0,d)$.  Let $\epsilon_o>0$ be a constant such that  $\ell_{c,d}\subset \ell_{c-2\epsilon_o ,d+2\epsilon_o }
 \subset\subset  \ell
$ and $\epsilon_o\leq \frac{d-c}{4}.$
Set
$$
E_{c-2\epsilon_o,d+2\epsilon_o}^{\delta_o}=[0,\delta_o]\times \ell_{c-2\epsilon_o ,d+2\epsilon_o}\subset \overline\Delta.
$$

We   assume that $E_{c-2\epsilon_o,d+2\epsilon_o}^{\delta_o}\bigcap (\partial \Delta\setminus \ell)=\emptyset.$

We use affine technique to prove the following proposition.
\begin{prop}\label{proposition_5.3.1}
There is a constant $\mff C_{8}>0$ independent of $k$ such that
$$d_{u\indexm}(p,q) \geq \mff C_{8}$$
for $k $ large enough. Here $d_{u\indexm}(p,q)$ denotes the geodesic distance from $p$ to $q$ with respect the Calabi metric $G_{u_k}$.
\end{prop}

  We introduce some notations. Let $u$ be a function of the sequence $u\indexm.$ Let $\Gamma$
 be a minimal geodesic from $(0,c)$ to (0,d) with respect to the Calabi metric $G_u$.
For any
$p^\circ\in \Gamma\setminus \ell$, denote $$d(p^\circ)=d(p^\circ,\partial\Delta).$$
Let $\check p\in \ell$ be the point such that
\begin{equation}\label{geqna_5.1}
\partial_2 u(\check p)=\partial_2 u(p^\circ).
\end{equation}

Let $L(p^\circ)$ be the geodesic arc-length  of the connected
component  containing $p^\circ$ of $\Gamma\cap B_{\tau d(p^\circ)}(p^\circ)$, where
$\tau$ is the constant in Lemma \ref{glemma_3}.  Then
$L(p^\circ)=2\tau d(p^\circ).$
 Denote  $$m(p^\circ)=\max_{q\in B_{\tau d(p^\circ)}(p^\circ)}
|\partial_2 u(q)-\partial_2 u(p^\circ)|,$$
$$\mathcal B(\check p )=\{p\in \ell \;|\; |\partial_2 u(p)-\partial_2
u(\check p )|< m(p^\circ)\}.$$

To prove Proposition \ref{proposition_5.3.1} we need the following lemma.
\begin{lemma}\label{glemma_4} Let $u\indexm \in \mc C^\infty(\Delta,v;K_o)$ be a sequence of
functions such that $u\indexm$ satisfies \eqref{eqnc_5.1}  and
$$\mc S(u\indexm)=K\indexm. $$
 Suppose that  $K\indexm$ uniformly  $C^3$-converges to  $K$ on $\overline\Delta$ and the geodesic  arc-length of $\Gamma_k$ converges to zero as $k\to \infty$. Then there is a positive constant  $C_2$ independent of k such that
 for any $u_k$ and any $p^\circ \in \Gamma_k\cap E^{\delta_o}_{c-\epsilon_o,d+\epsilon_o}, $
\begin{equation}
\int_{\mathcal B_k(\check p)} \sqrt{\partial ^2_{22}u_k} d\xi_2
\leq C_2 \tau d_k(p^\circ) \leq C_2 L_k (p^\circ),
\end{equation}
where $d_k(p^\circ)=d_{u\indexm}(p^\circ,\partial \Delta)$ and $\check p$ satisfies \eqref{geqna_5.1}.
\end{lemma}
{\bf Proof.} If the lemma is not true, there are a subsequence of points $p^\circ_k$
and a subsequence of functions  $u\indexm,$ still denoted by $p^\circ_k$ and $u\indexm $  to simplify  notations, such that
\begin{equation}\label{geqn_3.2}
\lim_{k\to\infty}\frac{\tau d_{k}(p^\circ_k) }{\int_{\mathcal B(\check p_k)} \sqrt{\partial ^2_{22}u\indexm} d\xi_2
}=0.
\end{equation}

Let $u$ be  a function of the sequence $u\indexm.$ Let  $\hat u=u-\nabla u(p^\circ)\cdot\xi+C$,   where $C$ is a constant such that $\inf_{\partial \Delta}\hat u=0. $
Then, $\hat  u(\check p) =\inf\limits_{\ell} \hat u, \hat u(p^\circ)=\inf\limits_{\Delta} \hat u.$  We claim that
\begin{equation}
\inf_{\partial \Delta}\hat u=\inf_{\ell} \hat u=\hat u(\check p),\;\;\;\;
\end{equation}
as $k$ is large enough.
\v\n
{\em Proof of the Claim.} If the Claim is not true, there are a subsequence of $\hat u\indexm$ and a sequence of points $\check q_k\in \partial \Delta\setminus \ell ,$ still denoted by $\hat u\indexm$ and $\check q_k,$ such that
$$\hat  u\indexm (\check q_k)=\inf_{\partial \Delta}\hat u\indexm.$$
Let $\alpha_k=-\inf\limits_{ \Delta}\hat u_k.$ Since  the geodesic  arc-length of $\Gamma_k$ converges to zero as $k\to \infty$ we have $\lim\limits_{k\to\infty}d_k( p^\circ_k )=0$. Then  by the interior regularity and choosing a subsequence, we have
\begin{equation}\label{eqnd_5.1}\lim_{k\to\infty}    p^\circ_k=\check p_\infty\in \ell_{c-\epsilon_o,d+\epsilon_o}.\end{equation} In fact, if $\check p_\infty\in \Delta,$ then $u_{k}$ $C^{3,\alpha}$ converges to a strictly convex function $u_{\infty}$ in the neighborhood of $\check p_\infty,$ and $d_{k}(p_{k}^{\circ})\geq \epsilon_{1}$ for some $\epsilon_{1}>0$ independent of $k$. It contradicts $\lim\limits_{k\to\infty}d_k( p^\circ_k )=0$.  It follows from \eqref{eqnd_5.1} and  Corollary \ref{corollary_uniform_control} that
\begin{equation}\label{eqnd_5.2}
\lim_{k\to\infty}     |u_k(p^\circ_k)-u_k(\check p_\infty)| =0.\end{equation}
By the convexity of $u$ and \eqref{eqnc_5.1} we have
$|\partial_2 u(p^\circ)|\leq C_1$
for some positive constant $C_1$ independent of $k$.
Then
\begin{eqnarray} &&\hat u (\check p_\infty )-\hat u (p^\circ)\nonumber \\ &=& u (\check p_\infty )-u (p^\circ) +
\partial_1  u (p^\circ ) \xi_1(p^\circ)-\partial_2 u (p^\circ)(\xi_2(\check p_\infty)-\xi_2 (p^\circ))\nonumber \\
\label{eqnd_5.3}
&\leq & u (\check p_\infty)-u (p^\circ)+C_1\left| (\xi_2(\check p_\infty )-\xi_2 (p^\circ))\right| +C\xi_{1}(p^{\circ}),
\end{eqnarray} as $k$ is large enough, where we used the fact $\partial_1 u (p ^\circ) <C , 0<\xi_1(p^\circ)<$diam$ (\Delta)$.
Combining \eqref{eqnd_5.1}, \eqref{eqnd_5.2}, \eqref{eqnd_5.3}, $\lim\limits_{k\to\infty}\xi_1(p_{k}^\circ)=0 $ and $ \hat u (\check p_\infty)> \hat u(\check q)>\hat u (p^\circ)$ we conclude that
  \begin{equation}\label{eqn_5.9a}
 \lim_{k\to\infty}\alpha_k=0.
\end{equation}

Let $\gamma$ be the line segment  connecting  $p^\circ$ and $\check q$. By the convexity we have for any $p\in \gamma\;$   \begin{equation}\label{geqn_3.4}-\alpha\leq \hat u(p)\leq 0 ,\;\;\;\;\; l(p)\leq u(p)\leq l(p)+\alpha,\end{equation}
 where  $l(p) =u(p^\circ)+\nabla u(p^\circ)\cdot (p-p^\circ).$
By choosing a subsequence  we can assume that  $\gamma_k$ converges to a line segment $\gamma_\infty.$  We can see that $\gamma_\infty\subset \ell,$
  otherwise, as $\lim_{k\to\infty}\alpha_k=0,$   $u_\infty$ is a linear function on $\gamma_\infty,$ it contradicts the strict convexity of $u_\infty$.
For any $p\in \gamma_\infty\cap \ell_{c-\epsilon_o,d+\epsilon_o}$ and $p_k\in \gamma_k$ with $\lim\limits_{k\to\infty}p_k=p,$ by  Corollary \ref{corollary_uniform_control} and the same calculation as \eqref{eqnd_5.3} we have
\begin{equation}\label{eqnd_5.4}0\leq \lim_{k\to\infty}(\hat u_k ( p)-\hat u_k (p_k))\leq \lim_{k\to\infty}(u_k ( p )-u_k (p_k))=0.
\end{equation} Let $\hat h=\hat u|_{\ell}.$
By \eqref{eqn_5.9a}, \eqref{geqn_3.4} and  \eqref{eqnd_5.4} we conclude that
\begin{equation}\label{eqnd_5.5}\lim_{k\to \infty} \hat h_k(p) =0,\;\;\; \forall p\in \ell_{c-\epsilon_o,d+\epsilon_o}\cap \gamma_{\infty}.
\end{equation}
On the other hand, since  $\hat h =h -\partial_2 u(p^\circ)\xi_2+C ,$ from Lemma \ref{lemma_5.2.1}
$\hat h_k$ converges to a strictly convex function $\hat h_{\infty}.$ It contradicts  \eqref{eqnd_5.5}.
The Claim is proved.

\v By a coordinate translation, we assume that $\xi(\check p)=0.$
  Consider the following affine transformation $T$,
$$\tilde \xi_1=a_{11}\xi_1,\; \tilde \xi_2=a_{21}\xi_1+a_{22}\xi_2,\;\;\;\;\;
\tilde u(\tilde \xi)=\lambda \hat u\left(\frac{\tilde \xi_1}{a_{11}},\frac{\tilde \xi_2}{a_{22}}
-\frac{ a_{21}\tilde\xi_1}{a_{11}a_{22}}\right).
$$
We choose $\lambda=[d(p^\circ)]^{-2} $ and
$$a_{11}=\frac{\sqrt{\lambda\det(\partial^2_{ij}\hat u)}}{\sqrt{\partial^2_{22}\hat u}}(p^\circ),\;\;
a_{21}=\frac{\sqrt{\lambda}\partial^2_{21}\hat u}{\sqrt{\partial^2_{22} \hat u}}(p^\circ),\;\;
a_{22}=\sqrt{\lambda \partial^2_{22} \hat u }(p^\circ).
$$
Denote by $\tilde p^\circ, \tilde {\check {p}},\cdots$  the image of $ p^\circ,  {\check {p}},\cdots$ under the affine transformation $T$.
Then by a direct calculation we have
$ \tilde \partial^2_{ij} \tilde u (\tilde p^\circ)=\delta_{ij},$
and for any $p,q$
\begin{equation}\label{geqn_3.5}
 \tilde\partial_2 \tilde u(\tilde p ) =
\frac{\lambda}{a_{22}} \partial_2  \hat u( p), \;\; \tilde\partial^2_{22} \tilde u(\tilde p ) =
\frac{\lambda}{a_{22}^2} \partial^2_{22}  \hat u( p),\;\; \tilde d_{\tilde u}(\tilde p ,\tilde q)=\sqrt{\lambda}d_{\hat u}(p,q) ,
\end{equation}
\begin{equation}\lim_{k\to\infty}\max |\mathcal S(\tilde u\indexm)|=\lim_{k\to\infty}\max \frac{ \mc S(u_k)}{\lambda_k}=\lim_{k\to\infty}\max d(p^\circ_k) \mc S(u_k)=0,
\end{equation}
where we denote $\tilde\partial_i \tilde u= \frac{\partial \tilde u}{\partial \tilde  \xi_i}, $ $\tilde\partial^{2}_{ij}\tilde u=\frac{\partial^2 \tilde u}{\partial \tilde \xi_i \partial \tilde \xi_j}$ and use $\mc S(u)=\mc S(\hat u)$.
In particular, $\tilde d(\tilde p^\circ)=1,$
and
\begin{equation}\label{geqn_3.7}
\lim_{k\to\infty}\frac{\tau \tilde d(\tilde p_k^\circ) }{\int_{\mathcal{\tilde B}(\tilde {\check p} _k)} \sqrt{\partial ^2_{22}\tilde u\indexm} d\tilde\xi_2
}=\lim_{k\to\infty}\frac{\tau  d(\tilde p_k^\circ)}{\int_{\mathcal{   B}( \check  p_k)} \sqrt{\partial ^2_{22}  u\indexm} d \xi_2
}=0,
\end{equation}
where $$
\tilde m(\tilde p^\circ)=\max_{B_{\tau}(\tilde p^\circ)}|\tilde\partial_2\tilde u |
,\;\;\;\;  \mathcal{\tilde B}(\tilde{\check p}) =\{q\in\t_2^\ast| |\tilde\partial_2 \tilde u|(q) <
\tilde m(\tilde p^\circ)
\}.
$$
On the other hand,   using  Lemma \ref{glemma_3} for $(\tilde u_k, \tilde p_k^\circ, \tilde{ \check p}_k ),$
 we conclude  that $\tilde u\indexm$ locally uniformly converges to $\tilde u_\infty$ and
\begin{eqnarray}
\tilde m(\tilde p^\circ) \leq C_1\inv,\;\;\;
\mathcal{\tilde B}(\tilde {\check p}) \subset S_{\tilde h}(\tilde {\check p},1) ,\;\;\;
|\tilde \xi_2(q)|\leq C_1/2,\;\;\forall\;\; q\in \mathcal{\tilde B}(\tilde p^\circ).
\end{eqnarray}
 Then
 $$\int_{\mathcal{\tilde B}(\tilde {\check p})} \sqrt{\tilde \partial ^2_{22}\tilde u\indexm} d\tilde \xi_2\leq
 \left(\int_{\mathcal{\tilde B}(\tilde {\check p} )} \tilde\partial ^2_{22}\tilde u\indexm d\tilde\xi_2   \right)^{\half}
 \left(\int_{\mathcal{\tilde B}(\tilde {\check p})}  d\tilde\xi_2\right)^{\half}\leq  \sqrt{2\tilde m(\tilde p^\circ)C_1} \leq 2.
 $$
It contradicts to \eqref{geqn_3.7}.  $\Box.$ \\

\noindent{\bf Proof of Proposition \ref{proposition_5.3.1}.} If the Proposition is not true, by choosing a subsequence we can assume that
\begin{equation}
\label{eqnd_5.7} \lim_{k\to\infty} L(\Gamma_k)=0,
\end{equation} where $L(\Gamma_k)$ denotes the geodesic arc-length of $\Gamma_k.$
 Moreover, we can assume that the Euclid measure of $\Gamma_k\cap \ell$ goes to zero as $k\to \infty.$ In fact, if the Euclid measure of $\Gamma_k\cap \ell$ has uniform positive lower bound, we can  get a contradiction easily  from $\partial^2_{22}u|_{\ell_{c-\epsilon_o,d+\epsilon_o}}\geq C.$

 There is a open set $U\subset \ell $  which contains  $\Gamma\cap\ell_{c-\epsilon_o,d+\epsilon_o}$, such that
 $\ell_{c-\epsilon_o,d+\epsilon_o}\setminus U$ is a compact set and the Euclidean measure of $U$  less than $\epsilon_o/2,$ as $k$ is large enough.
For any $\delta_{1}>0$ and $\epsilon_{o}>0,$ denotes
$$
L^{\delta_{1}}_{c-\epsilon_o}=\{(\xi_1, c-\epsilon_o)|0\leq \xi_1\leq \delta_1\},\;\;\;  L^{\delta_{1}}_{d+\epsilon_o}=\{(\xi_1,
d+\epsilon_o)|0\leq \xi_1\leq \delta_1\}.
$$
By Lemma \ref{lemma_grad_control}, there exists  a constant  $\delta_1>0$ independent of $k$  such that
\begin{equation}\label{eqnd_5.6}\max_{L^{\delta_{1}}_{c-2\epsilon_o}}\partial_2  u  \leq \min_{L^{\delta_{1}}_{c-\epsilon_o}} \partial_2  u  ,\;\;  \min_{L^{\delta_{1}}_{d+2\epsilon_o}} \partial_2  u  \geq \max_{L^{\delta_{1}}_{d+\epsilon_o}} \partial_2  u  .\end{equation}
Denote $p_1=(0,c-\epsilon_o),q_1=(0,d+\epsilon_o).$

  Since the geodesic  arc-length of $\Gamma_k$ converges to zero as $k\to \infty$,
by the interior regularity we can assume that $\Gamma \cap \{\xi\in E_{c-2\epsilon_o,d+2\epsilon_o}^{\delta_1}|\xi_1=\delta_1\}=\emptyset.$

We discuss three cases: \\
{\bf Case 1. }  $\Gamma \cap L^{\delta_{1}}_{c-2\epsilon}\neq \emptyset.$
Since $\partial_2 u$ is continuous on  $\Gamma \cap E_{c-2\epsilon_o,d+2\epsilon_o}^{\delta_1},$ by \eqref{eqnd_5.6}
we have
$$[\partial_2 u(p_1),\partial_2 u(p)]\subset \bigcup_{\xi\in \Gamma \cap E_{c-2\epsilon_o,d+2\epsilon_o}^{\delta_1}}\partial_2 u(\xi). $$
Hence
 $$\ell_{c-\epsilon_o,c}\setminus U\subset {\bigcup _{p^\circ\in \Gamma\setminus \ell}\mathcal B(\check p)},\;\;\;
 |\ell_{c-\epsilon_o,c} \setminus U|\geq  {\epsilon_o}/{2}
 .$$
{\bf Case 2. }  $\Gamma \cap L^{\delta_{1}}_{d+2\epsilon}\neq \emptyset.$
By the same argument of Case 1, we have
$$\ell_{d,d+\epsilon_o}\setminus U\subset {\bigcup _{p^\circ\in \Gamma\setminus \ell}\mathcal B(\check p)},\;\;\;
|\ell_{d,d+\epsilon_o} \setminus U|\geq  {\epsilon_o}/{2}
.$$
{\bf Case 3. }   $\Gamma\subset E_{c-2\epsilon_o,d+2\epsilon_o}^{\delta_1}.$
Since $\partial_2 u$ is continuous on $\Gamma,$ we have $$[\partial_2 u(p),\partial_2 u(q)]\subset \bigcup_{\xi\in \Gamma \cap E_{c-2\epsilon_o,d+2\epsilon_o}^{\delta_1}}\partial_2 u(\xi). $$
    Hence $$ \ell_{c,d} \setminus U\subset {\bigcup _{p^\circ\in \Gamma\setminus \ell}\mathcal B(\check p)},\;\;\;\; |\ell_{c,d} \setminus U|\geq  {\epsilon_o}/{2} $$  as $k$ is large enough.

    We prove the Case 3 (the proof of the other cases is the same).
     There are  finitely many  points $\{ p^\circ_l\}_1^N\subset \Gamma$ such that $\{\mathcal B(\check p_l)\}_1^N$ covers $\ell_{c ,d} \setminus U$ and
 $$\mathcal B(\check p_i)\cap\mathcal B(\check p_j)\cap \mathcal B(\check p_l) =\emptyset,$$
 for any different $i,j,l.$ Then
$$  B_{\tau d_i}( p^\circ_i)\cap   B_{\tau d_j}( p^\circ_j)\cap    B_{\tau d_l}( p^\circ_l) =\emptyset,$$
where $d_i=d( p^\circ_i,\partial \Delta).$ In fact,  if there exists a point $p^\ast\in   B_{\tau d_i}( p^\circ_i)\cap   B_{\tau d_j}( p^\circ_j)\cap    B_{\tau d_l}( p^\circ_l) $, let $q^\ast \in \ell$
such that $\partial_2 u(q^\ast)=\partial_2 u(p^\ast).$ By definition of  $\mathcal B(\check p_l)$
we have $q^\ast\in \mathcal B(\check p_i)\cap \mathcal B(\check p_j)\cap\mathcal B(\check p_l).$ We get a contradiction.
Therefore $L(\Gamma)\geq \half \sum\limits_{i=1}^N L(p^\circ_i)$.
By  Lemma  \ref{lemma_5.2.1} and Lemma \ref{glemma_4} we have
 $$L(\Gamma)\geq\sum_{i=1}^N \tau d_i \geq  \frac{1}{C_2}\int_{\ell_{c,d}\setminus U}\sqrt{\partial^2_{22} u}d\xi_2\geq \frac{\sqrt{\mff C_7}\epsilon_o}{C_2} .$$
  It contradicts  \eqref{eqnd_5.7}.
 We finish the proof of Proposition \ref{proposition_5.3.1}.  $\Box$

\v % Denote
%$$
%L^{\delta'}_c=\{(\xi_1, c)|0\leq \xi_1\leq \delta'\},\;\;\; L^{\delta'}_d=\{(\xi_1,
%d)|0\leq \xi_1\leq \delta'\}.
%$$
By the same argument, we have
\begin{prop}\label{proposition_5.3.2}
For any $\delta>0,$ there exists a constant $\mff C_{9}>0$ independent of $k$ such that
$$
d_{u_k}(p,q)>\mff C_{9}
$$
for $p\in L^{\delta}_c$ and $q\in L^{\delta}_d$.
\end{prop}

By the interior regularities and the same argument of Proposition \ref{proposition_5.3.1}, we conclude that
\begin{theorem}\label{theorem_5.3.3}
Let $p\in \ell_{c,d}.$ Suppose that  $D_a(p)\cap \partial \Delta\subset  \ell_{c-\epsilon_o,d+\epsilon_o}$ for some $0<a<\epsilon_{o}.$
Then there exists a small constant
$\epsilon$, independent of $k$, such that the intersection of the geodesic ball
$B\indexn_\epsilon(p)$ and $\Delta$ is contained in a  Euclidean   half-disk $D_a(p)\cap \Delta$.
\end{theorem}

%%%%%%%%%%%%%%%%%%%%%sect_6

\section{Upper bound of $H$}\label{sect_6}
Recall that $H=\frac{\det(g_{i\bar j})}{\det(f_{i\bar j})}.$ The following  theorem has been proved in \cite{CLS2}.
\begin{theorem}\label{theorem_6.1}  Let $(M,G)$ be a compact complex
manifold of dimension n with  K$\ddot{a}$hler metric $G$. Let $\omega_o$ be its K$\ddot{a}$hler
form. Denote
$$\mc C^{\infty}(M,\omega_o) = \{\phi \in C^{\infty}(M)| \omega_\phi = \omega_o +
\frac{\sqrt{-1}}{2\pi}\partial\bar{\partial}\phi >0\}.$$ Then for any $
\phi\in \mc C^{\infty}(M,\omega_o)$,  we have
\begin{equation}\label{eqn_6.1}
H \leq \left(2+\frac{\max\limits_{M}|\mc S(f)|}{n\mc {\dot {K}}}\right)^n\exp
\left\{2\mc {\dot {K}}(\max_M\{\phi\}-\min_M\{\phi\})\right\}.
\end{equation}
where $f=g+\phi$ and  $\dot{\mc K}= \max\limits_{M}\|Ric(g)\|^2_{g}, $ $Ric(g)$ denotes  Ricci tensors of the metric $\omega_g.$
\end{theorem}

\def \minequiv{\operatorname*{\sim}\limits^{min}}

\section{Lower bound of $H$}\label{sect_7}
In this section we will prove the following theorem.
\begin{theorem}\label{theorem_7.0.1}
Let $\Delta\subset \mathbb R^2$ be a Delzant polytope and
 $(M,\omega_o)$ be the associated  compact  toric surface. Let $K\in C^{\infty}(\bar\Delta)$ be an edge-nonvanishing
function and $u\indexm=v+\psi\indexm\in \mc
C^\infty(\Delta,v)$ be a sequence of functions with $\mathcal S(u\indexm)=K\indexm$. Suppose that
\begin{enumerate}\item[(1)] $K\indexm$ $C^{3}$-converges to $K$ on $\bar\Delta$;
\item[(2)]  $\max_{\bar \Delta}|u_k|  \leq \mff C_1$ ,\end{enumerate}
where $\mff C_1$ is a constant independent of $k$.
 Then there exists a constant $\mff C_{10}>1$ independent of $k$ such that for any
$k$
\begin{equation}\label{eqn_7.1}
 \mff C_{10}\inv\leq H_{f_k}\leq \mff C_{10}.\end{equation}
\end{theorem}

\v
The upper bound is proved.
Let $p\indexm\in\bar\Delta$  be the minimum point of
$H_{f\indexm}$, that is, for any $z_k\in \tau\inv_{f_k}(p_k),$
$H_{f_k}(z_k)=\min\limits_{M} H_{f_k}$. Let $p_\infty$   be the limit of $p\indexm$ (if necessary, by taking a subsequence to get the limit). Then by the
interior regularity,% theorem (cf. Theorem \ref{theorem_2.2.3}),
 we can assume that
$p_\infty \in \partial \Delta$.

\subsection{A subharmonic function}\label{sect_7.1}

Let $\cplane^2_\vartheta$ be a coordinate chart associated to the
vertex $\vartheta$. Let $\mc Q_\vartheta$ be the space of
coordinates of radius of $\cplane^2_{\vartheta}$. It is the first
quadrant of $\real^2$. We omit the index $\vartheta$ if
there is  no danger of confusion. We have a natural map
$$
\rho: \cplane^2\to \mc Q,\;\;\;
\rho(z_1,z_2)=(r_1,r_2)=(|z_1|,|z_2|).
$$
Since we consider  $\bb T^2$-invariant objects, we identify
$\cplane^2$ as $\mc Q$ in the following sense: when we
write a set $\Omega\subset \mc Q$, we mean $\rho\inv( \Omega)$.
We
also note that $\mc Q^\circ$ (the interior of $\mc Q$) is
identified with $\t$ in a canonical way: $x_i=2\log r_i$.

Introduce the notations in $\mc Q$
$$
\bx(a;b)=\{(r_1,r_2)|r_1\leq a, r_2\leq b\}.
$$
$$
\bx(a_1,a_2;b_1,b_2)=\{a_1\leq r_1\leq a_2, b_1\leq r_2\leq b_2\}.
$$
Let $\mc B_\vartheta=\bx(1,1)$ in $\mc Q_\vartheta$. Its boundary
consists of two parts,
 $E^i_\vartheta=\{(r_1,r_2)| |r_i|=1, |r_{3-i}|\leq 1 \}, i=1,2$.
(Here, by the boundary we mean the boundary of $\rho\inv(\mc B)$
in $\cplane^2$. Hence the boundaries of the box located on the
axis are in fact the interior of the complex manifold.)

For a toric surface, we have the following simple lemma.
\begin{lemma}\label{lemma_7.1.1}
(i) All $\mc B_\vartheta$'s in $M$ intersect at one point, i.e,
$(1,1)$ in each $\mc B_\vartheta$; (ii) For any two vertices
$\vartheta$ and $\vartheta^o$  next to each other, $\mc B_\vartheta$
and $\mc B_{\vartheta^o}$ share a common boundary; (iii)
$M=\bigcup_{\vartheta}  \mc B_\vartheta .$
\end{lemma}
{\bf Proof. }Let $\vartheta$ and $\vartheta^o$ be two vertices
that are next to each other. Let $\ell$ be the edge connecting
them. We put $\Delta$ in the first quadrant of $\t^\ast$ as the
following: {\em (1) $\vartheta$ at the origin; (2) $\ell$ on the
$\xi_2$-axis; (3) the other edge $\ell_\ast$ of $\vartheta$ on the
$\xi_1$-axis; (4) suppose that $\vartheta^o=(0,c_o)$ and its other
edge $\ell^\ast$  is given by the equation $ \xi_2=a_o\xi_1+c_o $
for some integer $a_o$.}

The edges $\ell_\ast$ and $\ell$  is of the
form $\vartheta+te_1$ and $\vartheta+te_2,t\in \mathbb R$ respectively. Here $\{e_1,e_2\}$ is a basis of $\mathbb Z^2$.
Similarly,
$\ell^\ast$ and $\ell$ is of the form $\vartheta^o+te_1^o$ and $\vartheta^o+te_2^o,t\in \mathbb R$ respectively, and $\{e_1^o,e_2^o\}$ is a basis of $\mathbb Z^2.$
Then $$e_1^o=e_1+a_oe_2,\;\;\; e_2^o=-e_2.$$
For any point $p\in \Delta,$
we have
$$p=(e_1,e_2)(\xi_1,\xi_2)^t=(e_1^o,e_2^o)(\xi_1^o,\xi_2^o)^t+(e_1,e_2)(0,c_o)^t,$$
where $A^t$ denotes the transpose of a matrix $A$.
Hence the coordinate transformation
between $(\xi_1,\xi_2)$ and $(\xi_1^o,\xi_2^o)$ is
$$\xi_1^o=\xi_1,\;\;\;\;\xi_2^o=a_o\xi_1-\xi_2+c_o,$$
and the coordinate transformation
between $(x_1,x_2)$ and $(x_1^o,x_2^o)$ is
$$x_1^o=x_1+a_ox_2,\;\;\;\;x_2^o=-x_2.$$
Then the coordinate transformation
between $\cplane^2_\vartheta$ and $\cplane^2_{\vartheta^o}$ is
given by
\begin{equation}\label{eqn_7.2}
z_1^o= z_1z_2^{a_o}, z_2^o=z_2\inv.
\end{equation}
By this, we find that $\mc B_\vartheta$ and $\mc B_{\vartheta^o}$
intersect at the common boundary $E_\vartheta^2=E_{\vartheta^o}^2$.
The rest of the  facts of the lemma can be derived easily as well. q.e.d.
\v Let $\mc E$ be the collection of all $E^i_\vartheta$'s.

Now consider an element $\mathsf{f}\in \mc C^\infty(M,\omega;
K_o)$. Let $f_\vartheta$ be its restriction to  $\CHART_\vartheta$.
We introduce a {\em subharmonic} function
$$
F_\vartheta=\log W_\vartheta+Nf_\vartheta.
$$
\begin{lemma}\label{lemma_7.1.2}
If $N\geq \max|\mc S(f_\vartheta)|+1$, $\square F_\vartheta>0$.
Hence the maximum  of $F_\vartheta$ on $ \mc B_\vartheta $ is  achieved on
$ E^1_\vartheta\cup E^2_\vartheta $. Here $\square$ denotes the complex Laplacian operator of the metric $\omega_{f}.$
\end{lemma}
{\bf Proof. }By a direct computation,
$$
\square F_\vartheta=-\mc S(f_\vartheta) +2N>0.
$$
q.e.d.

\v
\begin{lemma}\label{lemma_7.1.3}
All $F_\vartheta$ on $ \mc B_\vartheta $ form a continuous function
$\mathsf F$ on $M$.
\end{lemma}
{\bf Proof. }Let $\vartheta$ and $\vartheta^o$ be two vertices
that  are next to each other as in the proof of Lemma \ref{lemma_7.1.1}.
 By
a direct calculation, we have
\begin{eqnarray}
f_{\vartheta^o}-f_{\vartheta}&=&\sum  ( x_i^o   \xi_i^o-u)- \sum (x_i \xi_i-u)  \nonumber\\
&=& \sum x_i^o \xi^o_i - (x_1^o+a_o x^o_2)\xi^o_1 -(- x^o_2)(a_o\xi^o_1- \xi^o_2+c_o) \nonumber\\
&=&\label{eqnb_8.1} c_o \log |r^o_2|^2.\\
\label{eqnb_8.2}
\log W_{\vartheta^o}-\log W_{\vartheta}&=& \log \left|\frac{\partial z_i}{\partial z^o_j}\right|^2 =(a_o-2)\log|r_2^o|^2.
\end{eqnarray}
From this we conclude that $ F_{\vartheta^o}$ and $F_{\vartheta}$
match on their common boundary (where $r_2=r_2^o=1$). Hence,  all $F_\vartheta$'s form a function on $M$. q.e.d \v
Hence $\mathsf F$ is a continuous function on $M$ and piecewise strict
subharmonic.
%\begin{defn}\label{definition_7.1.4}
%A point $p\in E\subset \mc E$ is called a local maximum point of
%$\mathsf F$ if it is a maximum in $ \mc B_\vartheta \cup \mc
%B_{\vartheta^o} $, where the vertices $\vartheta,\vartheta^o$ are
%chosen such that $\mc B_\vartheta\cap \mc B_{\vartheta^o}=E$.
%\end{defn}
%\begin{corollary}\label{corollary_7.1.5}
%The set of local maximum points of $\mathsf F$ is nonempty.
%\end{corollary}

\def \fkz{\mathfrak{z}}
\def \sff{\mathsf{f}}
 Recall that $W_{g_\vartheta}=\det((g_\vartheta)_{i\bar j}).$
\begin{lemma}\label{lemma_7.1.6}
Let $q_o$ be the maximum point of $\mathsf F.$ Then there is a constant $\mathsf C_{11}>0$ independent of $k$ such that
 $H(q_{o})\leq \mathsf C_{11} \min_{M} H$
\end{lemma}
{\bf Proof. }Let $z_o$ be the minimal point of $H$. Suppose that
it  is in $\mc B_{\vartheta_0}$ for some vertex $\vartheta_0$.
Suppose that $\mathsf F$ achieves its maximum at $q_{o}$ in  $\mc B_{\vartheta}$. %We claim that $p_0\minequiv z_o$: in fact,
 by the assumption, $
F_{\vartheta}(q_{o})\geq F_{\vartheta_0}(z_o). $ Explicitly, this is
 $$\log W_{\vartheta}(q_{o})+Nf_{\vartheta}(q_o)\geq \log W_{\vartheta_{0}}(z_o)
+Nf_{\vartheta_{0}}(z_o).$$ Hence, by the definition of $H$, we have
\begin{eqnarray*}
&&\log H(q_o)-\log W_{g_{\vartheta}}(q_o)-N(f-g)_{\vartheta}(q_o) -Ng_{\vartheta}(q_o) \\
&\leq& \log H(z_o) -\log
W_{g_{\vartheta_0}}(z_o)-N(f -g)_{\vartheta_0}(z_o)-Ng_{\vartheta_0}(z_o).
\end{eqnarray*}
In $\mc B_{\vartheta}$ (resp.$\mc B_{\vartheta_0}$), $|g_{\vartheta}|$ (resp.$|g_{\vartheta_0}|$)  and $|\log W_{g_{\vartheta}}|$ (resp.$|\log W_{g_{\vartheta_0}}|$) are uniformly
bounded. Note that $\|f_\bullet-g_\bullet\|_{L^{\infty}}=\|u-v\|_{L^{\infty}(\Delta)}\leq \mathcal C_{1}.$ Therefore, there exists a constant $C$ such that
$$\log H(q_o)\leq\log H(z_o)+C .$$ This implies the claim. \;\;\; q.e.d.

\v
Let
$$
A^\delta_\vartheta= \bx(\delta,\delta\inv;\delta,\delta\inv)
\subset \mc Q_\vartheta.$$ Set $ M^\delta\subset M$ to
 be the union of all $ A^\delta_\vartheta.$ For simplicity we choose $\delta=\frac{1}{1000}.$
By the interior regularity we can assume that $q_{o}\notin M^{\delta}.$

\subsection{Proof of Theorem \ref{theorem_7.0.1}}
\label{sect_7.3}
To prove Theorem \ref{theorem_7.0.1} we need the
following lemma (cf. Lemma 7.16 in \cite{CLS}).
\begin{lemma}\label{lemmaa_7.5.1}
Let $z^\ast\in Z$. Let $f$ be a function in the sequence $f\indexm$. Suppose that in $B_{2a}(z^\ast)$, $ \mc K\leq
C_0$ for some constant $C_0>0$ independent of $k.$ Then there is a
constant $c>0$, independent of $k$, such that there is a point $z^o$
in $B_a(z^\ast)$ satisfying
$$
d(z^o, B_{2a}(z^\ast)\cap Z)=c.
$$
Obviously $c\leq a$. Hence $d(z^o,Z)=c$.
\end{lemma}

\v\n {\bf Proof of Theorem \ref{theorem_7.0.1}.}  The upper bound is proved in Theorem \ref{theorem_6.1}.   Let $f$ be a function in the sequence $f\indexm$ that satisfies
the conditions of Theorem \ref{theorem_7.0.1}.
Let $q_o$ be the point in Lemma \ref{lemma_7.1.6} and $p_{o}=\tau_{f}(q_{o}).$

Recall that the vertices and  the edges of $\Delta$ are denoted by
$$
\{\vartheta_0,\ldots,\vartheta_d=\vartheta_0\}, \;
\{\ell_0,\ell_1,\ldots, \ell_{d-1},\ell_d=\ell_0\},\;\;\mbox{where }\vartheta_{i}=\ell_{i}\cap \ell_{i+1},\; 0\leq i\leq d.
$$
 For any edge $\ell_{i},0\leq i\leq d,$ let  $\xi^{(\ell_{i})}\in\ell_{i} $
and
$
\mc D^{\ell_{i}}:= D_{ \epsilon}(\xi^{(\ell_{i})})\cap  \bar\Delta
$ be  a half $ \epsilon$-disk such that
\begin{equation}\label{jv_eqn7.5}
d_E(\mc D^\ell,\partial \Delta\setminus \ell)>\epsilon,\;\;\;\;\;\;|K|>\delta_o>0, \;\;\; \mbox{ on }\;\;\; \mc D^\ell.
\end{equation} for some $\delta_o>0$ independent of $k$.  Let $o$ be the mass center of $\Delta.$ For any $0\leq i\leq d,$ we denote by $o\xi^{(\ell_{i})}\vartheta_{i}\xi^{(\ell_{i+1})}$ the quadrilateral with the vertices $o$, $\xi^{(\ell_{i})}$, $\vartheta_{i}$ and $\xi^{(\ell_{i+1})}$. Then $$\bar \Delta =\bigcup_{i=1}^{d}o\xi^{(\ell_{i})}\vartheta_{i}\xi^{(\ell_{i+1})}.$$

Without loss of generality we can assume that $p_{o}\in o\xi^{(\ell_{1})}\vartheta_{1}\xi^{(\ell_{2})}$. Consider the coordinates $\mathbb C_{\vartheta_1}.$ Denote $\Omega_{\vartheta_1}=\tau_f^{-1}(o\xi^{(\ell_{1})}\vartheta_1\xi^{(\ell_{2})}).$ Let  $\fkz_{\ast}\in \partial \Omega_{\vartheta_1}$ be the point such that
$$F_{\vartheta_1}(\fkz_{\ast})=\max_{\Omega_{\vartheta_1}}F_{\vartheta_1}\geq F_{\vartheta_1}(q_{o}),$$ where  we used the fact that $F_{\vartheta_1}=\log W_{\vartheta_1} +Nf_{\vartheta_1}$ is a subharmonic function.   Note that $f_{\vartheta_1}$ and $\det((g_{\vartheta_1})_{k\bar l})$ are uniform bounded in $\Omega_{\vartheta_1}.$
It follows that
\begin{equation}\label{equiv_H}
H(\fkz_{\ast})\leq C_{1}'H(q_{o})\leq C_{1} \min H
\end{equation}
for some constants $C_{1}'>0$, $C_1>0$.
Denote $p_{\ast}=\tau_{f}(\fkz_{\ast}).$
Assume that
$d_{E}(p_{\ast},\xi^{(\ell_{1})} )\leq \frac{\epsilon}{4} $ or $d_{E}(p_{\ast},\xi^{(\ell_{2})} )\leq \frac{\epsilon}{4},$ otherwise, the theorem follows from the interior regularity and \eqref{equiv_H}. We assume  $d_{E}(p_{\ast},\xi^{(\ell_{1})} )\leq \frac{\epsilon}{4} $.

  By Theorem \ref{theorem_5.3.3}   there exists a constant $a>0$  independent of $k$ such that   $B_{2a}(\fkz_\ast)\subset \tau_{f}\inv(\mathcal D^{\ell_{1}})$.    Denote $W=\det(f_{i \bar j})$ and $W_{g}= \det(g_{i \bar j}).$   Since $W_{g}$ is uniform bounded in $\tau_{f}\inv(\mathcal D^{\ell_{1}}),$ it follows from \eqref{equiv_H} that
\begin{equation}\label{eqnc_7.1}
W(z)\leq    N_1 W(\fkz_\ast),\;\;\; \forall z\in B_{2a}(\fkz_\ast),
\end{equation}
 for some constants $ N_1>0$ independent of $k$.
  Applying Theorem \ref{theorem_3.4.2}, we can find a constant
$C_2>0$ independent of $k$ such that  in
 $
B_{a}(\fkz_\ast),$
 $$\left[\frac{W}{W(\fkz_\ast)}\right]^{\half}(\Psi+\mc K) \leq C_2.$$
  Notice that
$\Psi=\|\nabla \log W\|_f^2$.
 Let $a'=\min(a,\frac{1}{2\sqrt{C_2}})$.
 Then for any $z\in B_{a'}(\fkz_\ast),$
 \begin{equation}\label{eqn_7.6}
 \half\leq \left[\frac{W(z)}{W(\fkz_\ast)}\right]^{\frac{1}{4}}\leq\frac{3}{2},\;\;\;\;\;\mc K \leq 4C_2.
 \end{equation}
 On the other hand, by   Lemma \ref{lemmaa_7.5.1}, there is a
 $p' \in \tau_{f}(B_{a'/2}(\fkz_\ast))$ such that
 $d(p' ,\partial \Delta)= c'$ for some constant $c'>0$ independent of $k$.
  \v

{\em Claim: There is a constant $C_{3}>0$ independent of k, such that  $\xi_1(p' )\geq C_{3}.$}

\v\n{\em Proof of the Claim.} If the Claim is not true,  $\lim\limits_{k\to\infty}\xi_1(p'_k)=0. $ By \eqref{jv_eqn7.5}, the convexity of $u$ and $|\max\limits_{\bar \Delta}u-\min\limits_{\bar \Delta}u|\leq \mathcal C_{1}$ we have in $ \mathcal D^{\ell_{1}},$
$$|\partial_{2} u|(p')\leq \mathcal C_{1}\epsilon\inv,\;\;\;\; \partial_{1} u (p')\leq \mathcal C_{1}\epsilon\inv.$$
 Without loss of generality  we can assume that $\partial_2 u_k(p_k')=0$ (for general case, since $|\partial_2 u(p')|\leq \mathcal C_{1}\epsilon\inv$
 we can use $u-\partial_2 u(p')\xi_2$ instead of $u,$   and use the same argument).
 Consider the function $u^\ast=u-\partial_1u(p')\xi_1.$ Then $u^\ast(p')=\inf u^\ast.$
 By Lemma 7.5 in \cite{CLS} we have
  \begin{equation}\label{eqna_8.9} \inf_{\ell}u^\ast- u^\ast(p')\geq C_4>0
   \end{equation}
   for some constant $C_4>0$ independent of $k$. We discuss two cases. \\

   \n
     {\bf Case 1} $\partial_1u(p')<0.$ Then
   $$\inf_{\ell}u -u (p')  \geq  \inf_{\ell}u^\ast- u^\ast(p')\geq C_4.$$
    {\bf Case 2} $0\leq \partial_1u(p')<\mathcal C_{1}\epsilon\inv.$ Then by $\lim\limits_{k\to\infty}\xi_1(p'_k)=0$ we have
    $$ \inf_{\ell}u -u (p')  \geq  \inf_{\ell}u^\ast- u^\ast(p')-\frac{C_4}{2}\geq \frac{C_4}{2},$$
    as $k$ large enough.

  For two cases we have $ \inf_{\ell}u -u (p')\geq \frac{C_4}{2}.$
    By this and  Proposition \ref{proposition_5.1.1}, we
   get a contradiction. The Claim is proved.

\v Let $z' \in \tau_{f}^{-1}(p' )\in B_{a'}(\fkz_\ast)$ be the
corresponding point of $p' $.
 Following from the Claim and the  interior regularity,
$W(z' )$ is bounded as   $C\inv<W (z')<C$ for some
constant $C$ independent of $k$.   By \eqref{eqn_7.6}, $W (\fkz_{\ast  })$ is
bounded above, therefore $H(\fkz_{\ast })>C_5>0$ for some constant $C_5$ independent of $k$.  This
completes the proof of Theorem \ref{theorem_7.0.1}.
 $\blacksquare$

%%%%%%%%%%%%%%%%%%%sect_8

\section{Proof of Theorem \ref{theorem_1.1}}\label{sect_8}

\subsection{The continuity method}\label{sect_2.2}
We argue the existence of
the solution to \eqref{eqn_1.2} by the standard continuity method.

Let $K$ be the scalar function on $\bar\Delta$ in Conjecture
\ref{conjecture_2.1.3} and suppose that there exists a constant
$\lambda>0$ such that $\Delta$ is $(K,\lambda)$ stable.

Let $I=[0,1]$ be the unit interval. At $t=0$ we start with a known
metric, for example $\omega_o$ (cf. Remark \ref{remark_2.2.1}).
Let $K_0$ be its scalar curvature
on $\Delta$. Then $\Delta$ must be $(K_0, \lambda_0)$ stable for some constant
$\lambda_0>0$ (cf.\cite{CLS5}). At $t=1$, set
$(K_1,\lambda_1)=(K,\lambda)$. On $\Delta,$ set
$$
K_t=tK_1+(1-t)K_0, \;\;\; \lambda_t=t\lambda_1+(1-t)\lambda_0.
$$
It is easy to verify that $\Delta$ is $(K_t,\lambda_t)$ stable.
 Set
$$
\Lambda=\{t|\mathcal  S(u)= K_t \mbox { has a solution in }  \mathcal C^{\infty}(\Delta, v) .\}
$$
Then we should show that $\Lambda$ is open and closed. Openness is
standard by using Lebrun and Simanca's argument (cf.\cite{D4,L-Ss}).
 It
remains to get a priori estimates for solutions $u_t=v+\psi_t$ to
show  closedness.
However, for technical reasons, we are only  able to prove  closedness under the condition that $K$ is an edge-nonvanishing
function.

\begin{theorem}\label{theorem_2.2.2}
Let $\Delta\subset \mathbb R^2$ be a Delzant polytope and
 $(M,\omega_o)$ be the associated  compact  toric surface. Let $K\in C^{\infty}(\bar\Delta)$ be an edge-nonvanishing
function and $u\indexm=v+\psi\indexm\in \mc
C^\infty(\Delta,v)$ be a sequence of functions with $\mathcal S(u\indexm)=K\indexm$. Suppose that
\begin{enumerate}\item[(1)] $K\indexm$ $C^{3}$-converges to $K$  on $\bar\Delta$;
\item[(2)]  $
\max_{\bar \Delta}\left|u_k\right| \leq \mff C_1$ ,\end{enumerate}
where $\mff C_1>0$ is a constant independent of $k$.
 Then  there is a subsequence of  $\psi\indexm$ which  $C^{6,\alpha}$-converges to
a function $\psi\in  C^{6,\alpha}(\bar \Delta)$ with $\mathcal
S(v+\psi)=K$.
\end{theorem}
We can restate Theorem \ref{theorem_2.2.2} as follows:

\v\n {\bf Theorem \ref{theorem_2.2.2}'} {\em Let $\Delta\subset \mathbb R^2$ be a Delzant polytope and
 $(M,\omega_o)$ be the associated  compact  toric surface. Let $K\in C^{\infty}(\bar\Delta)$ be an edge-nonvanishing
function and $\phi\indexm\in
C^\infty_{\bb T^2}(M)$ be a sequence of functions with $$\mathcal S(f_k)=K\indexm\circ \nabla^{f\indexm},\;\;\; \omega_{f\indexm}>0 ,$$ where  $f\indexm=g+\phi\indexm  $ and $g$ is the potential function of $\omega_o$.      Suppose that
\begin{enumerate}\item[(1)] $K\indexm$ $C^{3}$-converges to $K$   on $\bar\Delta$;
\item[(2)]  $\max_{M}|\phi_k| \leq \mff C_1$,\end{enumerate}
where $\mff C_1>0$ is a constant independent of $k$.
 Then  there is a subsequence of  $\phi\indexm$ which   $C^{6,\alpha}$-converges to
a function $\phi\in  C^{6,\alpha}_{\mathbb T^2}(M)$ with $\mathcal
S(g+\phi)=K\circ\nabla^f$.
 }
\v
The estimates for $\psi\indexm$ will be established by the following steps:
\v\n
 {\em Interior estimates: }
Donaldson proved the interior regularity for the Abreu equation
when $n=2$. In \cite{CLS2}, we proved the interior regularity for
the Abreu equation  of toric manifolds for arbitrary $n$ assuming
the $C^{0}$ estimates.
%Due to Theorem \ref{theorem_2.1.8}, our result gives an alternative proof of the interior regularity for the  $n=2$ case. The statement is
%\begin{theorem}\label{theorem_2.2.3}
%Let $u\in \mc C^\infty(\Delta,v),k\in \mathbb Z_{\geq 0}$ and $\mc S(u)=K$. If $\|K\|_{C^{k}(\bar\Delta)}\leq C_o$ and$\max_{\bar \Delta}\left|u\right|\leq \mff C_1$, then for any $\Omega\subset\subset \Delta$, there is a constant $\mc C_1'$, depending on $\Delta$, $C_o,\mff C_1$ and $d_E(\Omega,\partial\Delta)$, such that $$\|u\|_{C^{k+3,\alpha}(\Omega)}\leq \mc C_1'.$$ \end{theorem} Here, $d_E$ is the Euclidean distance function.

\v\n
{\em Estimates on edges: }This is the most difficult part.
On each $\ell$, we use the point $\xi^{(\ell)}$ in the interior
of $\ell$. By the condition,
there exists a half $\epsilon$-disk
\begin{equation}\label{eqn_2.1}
\mc D^\ell:= D_\epsilon(\xi^{(\ell)})\cap  \bar\Delta,\;\;\;\;\;D_{2\epsilon}(\xi^{(\ell)})\cap \partial \Delta\subset \ell
\end{equation}
such that $ K$ is non-zero on this half-disk. Hence there exists
a constant $\delta_o>0$ such that
\begin{equation}\label{eqn_2.2}
|K|>\delta_o, \;\;\; \mbox{ on }\;\;\; \mc D^\ell.
\end{equation}
In \S\ref{sect_8.1} we will show the  regularity on
 a neighborhood of $\xi^{(\ell)}$ that lies  inside $\mc D^{\ell}$.

\v\n
{\em Estimates on vertices: }Once the first two steps are completed,
the regularity on a neighborhood of vertices is based on  a subharmonic function.
This is done in \S\ref{sect_8.2}.

\begin{remark}\label{remark_2.2.1}
Let $K$ be an edge-nonvanishing function on $\Delta$. Suppose that
it satisfies \eqref{eqn_2.2}.
By the computation in Proposition \ref{proposition_10.1.1}, we find that we may modify $\omega_o$ to a new
form $\tilde\omega_o$ (equivalently, from $g$ to $\tilde g$)
such that the scalar curvature $\tilde K_0$ also satisfies $\eqref{eqn_2.2}$ and $$\tilde K_0K>0, \;\;\; \mbox{ on }\;\;\; \mc D^\ell.$$
Hence we can assume that the whole path $K_t$ connecting $\tilde K_0$ and $K_1=K$
satisfies \eqref{eqn_2.2}.
\end{remark}

%\appendix
%\subsection{Construction of Scalar Curvature}\label{sect_10}

 Fix a point
$q_\ell$ inside $\ell$. We assign a sign $\text{sign}(\ell)=\pm 1$ to
$q_\ell$. Then we conclude that
\begin{prop}\label{proposition_10.1.1}
There is a potential function $v$  on $\Delta$ such that for any $\ell$
$$
{\rm sign}(\ell)K(q_\ell)>0,
$$
where $K=\mc S(v_{o})$.
\end{prop}

 Let $\ell$ be a edge. We choose a coordinate on $\t^\ast$ such that $\ell=\{\xi|\xi_1=0\}$
and $q_\ell$ is the origin.
For any $\delta>0,$ let
$$
\Omega^1_\delta(\ell)=:\{\xi\in \Delta| \xi_1\leq
\delta \},\;\;\;
\Omega_\delta^2(q_\ell):
=\{\xi\in \Delta| |\xi_2|\leq\frac{\delta}{2}\},
$$
and $\Omega_\delta(q_\ell)=\Omega^1_\delta(\ell)
\cap \Omega_\delta^2(q_\ell)$. Set $\Omega^c(\ell)=\Delta\setminus (\Omega^1_{2\delta}(\ell)
\cup \Omega^2_{2\delta}(\ell))$.
\v
 We need the following lemma.
\begin{lemma}\label{lemma_10.1.3}
There is a convex function $u_\ell$ such that   ${\rm sign}(\ell)\mc S(u_\ell)>0$   in $\Omega_\delta(q_\ell)$;
and
 $u_\ell$ is linear on $\Omega^c(\ell)$.
\end{lemma}
{\bf Proof. } Consider the convex function
$u=\xi_1\log\xi_1 +\psi$ where $\psi$ is a function of $\xi_2$.
We compute $\mc S(u)=-\sum u^{ij}_{ij}$.
By a direct
calculation, we have
\begin{equation}
\mc S(u)= -\left(\frac{1}{\psi_{22}}\right)_{22}.
\end{equation}
Set $\Psi=(\psi_{22})\inv$. We consider $\Psi$ to be a
function of the following form:
\begin{equation}
\Psi= a\xi_2^2+c.
\end{equation}
Then $\mc S(u)=-2a$. We may choose $a$ to have the right sign.
Now we have
$$
\phi:=\psi_{22}=\frac{1}{a\xi_2^2+c}.
 $$
By choosing $c$ large, $\psi_{22}$ is positive.

\v Now we construct two convex functions $\alpha(\xi_1)$ and $\beta(\xi_2),$ such that
\begin{itemize} \item $\alpha(\xi_1)=\xi_1\log\xi_1$ when $\xi_1\leq \delta$
and is a linear function when $\xi_1\geq 2\delta$;
\item
and  $\beta(\xi_2)=\psi$ when $|\xi_2|\leq \delta/2$
and is  a linear function when $|\xi_2|\geq \delta$, where
$\psi$ is the function as above.
\end{itemize}
We use a cut-off function to modify $\phi$ to be a non-negative function $\tilde \phi$
 such that $\tilde\phi(\xi_2)=\frac{1}{a\xi_2^2+c}$ when $|\xi_2|\leq \delta$
and vanishes when $\xi_2\geq 2\delta$. Then $\beta$ can be constructed from $\tilde\phi$
such that $\beta''=\tilde\phi$ and $\beta(0)=\psi(0),\beta'(0)=\psi'(0)$.
By the same method we can construct the function $\alpha(\xi_1).$\\
\quad \quad Let
  $u_\ell= \alpha(\xi_1)+\beta(\xi_2).$  Then  $u_{\ell}$ satisfies the lemma. q.e.d.
\v\n
{\bf Proof of Proposition \ref{proposition_10.1.1}.} It is not hard to choose $\delta$ and arrange the coordinate system when
we construct $\psi$ such that
\begin{equation}
\Omega_\delta(q_\ell)\subset \bigcap_{\ell'\not=\ell}
\Omega^c(\ell').
\end{equation}
Let
$
\tilde u=\sum_{\ell} u_{\ell},
$ where $u_{\ell}$ is the function in Lemma \ref{lemma_10.1.3}. Then on $\Omega_\delta(q_\ell)$,
$
\mc S(\tilde u)=\mc S(u_\ell).
$
Hence, in $\Omega_\delta(q_\ell)$,
$$
\textrm{sign}(q_\ell)\mc S(\tilde u)>0.
$$
However, $\tilde u$ is not strictly convex.
Let
$$
u=\tilde u+\epsilon(\xi_1^2+\xi_2^2).
$$
For $\epsilon>0 $  small enough,
$
\mc S(u)$ is a small perturbation of $\mc S(\tilde u)=\mc S(u_\ell)
$
near $q_\ell$. The proposition is then proved. q.e.d.

\begin{remark}\label{remarka_2.2.1}
%Since $\Lambda$ is an open set, there is  $t_0>0$ such that  $[0,t_0)\subset \Lambda. $ Obviously $[0,\frac{t_0}{2}]\subset \Lambda.$
 For any $t\in [0,1],$ $\Delta$ is $(K_t,\lambda')$-stable, where $\lambda'=\min\{\lambda_{0} ,\lambda_1\}.$
Let $u_t$ be a solution of the equation $S(u)=K_t.$ Applying Theorem \ref{theorem_2.1.8} we have
$$|\max\limits_{\bar \Delta} u_t-\min\limits_{\bar \Delta} u_t|\leq \mathcal C_1,\;\;\forall\; t\in [0,1],$$ where $\mathcal C_1$ is a constant depending only on  $\lambda'$, $\Delta$ and $\|K_0\|_{C^{0}}+\|K_1\|_{C^{0}}.$
\end{remark}

 It remains to prove the regularity on divisors.
\def \fkz{\mathfrak{z}}
\subsection{Regularity on edges}
\label{sect_8.1}
Let $\ell$ be any edge and $\xi^{(\ell)}\in \ell$ such that
$|\mc S(\xi^{(\ell)})| >0$.
Recall that $f_{k}=g +\phi_k,$ where $f_k,g$ are Legendre transform of $u_k,v$ respectively; and $\phi_k\in C^\infty_{\mathbb T^2}(M).$
By Remark \ref{remark_2.2.1}  we can assume
that $|K_k|>\delta >0$ in  $$\Omega:=\{(z_1,z_2)|\log|z_1|^2\leq \half ,|\log |z_2|^2|\leq 1\; \},\;\;D_{2a}(\xi^{(\ell)})\subset \tau_{f_k} (\Omega)$$
and  $|z_1(\fkz^{(\ell)}_k) |=0,|z_2(\fkz^{(\ell)}_k)|=1,$
where $\delta,a$ are positive  constants independent of $k,$
  $\fkz^{(\ell)}_k\in Z_\ell$  whose image of the moment map is $\xi^{(\ell)}$.

  We omit the index $k$ if
there is  no danger of confusion. By Theorem \ref{theorem_5.3.3}, we conclude that
there is a constant $\epsilon>0 $ that is independent of
$k$   such that
$B _\epsilon(\xi^{(\ell)})\cap \Delta\subset D_a(\xi^{(\ell)})\cap \Delta$. Then
$B _\epsilon(\fkz^{(\ell)})$ is uniformly bounded. Hence, on this domain, we
assume that all data of ${g_\ell}$
are bounded.

Note that $W_{g_\ell}$ is bounded on $B _\epsilon(\fkz^{(\ell)})$. By  Theorem \ref{theorem_7.0.1} we have on $B _\epsilon(\fkz^{(\ell)})$
\begin{equation}\label{jv_add8.1}
C_1\inv\leq W \leq C_1
\end{equation}

It follows from Lemma \ref{lemma_5.2.1} that
$$ \partial ^{2}_{22} h|_{D_a(\xi^{(\ell)})\cap \ell}\geq \mathsf C_{7}.$$
Then by Theorem \ref{theorem_3.4.2} and \eqref{jv_add8.1} we conclude that there is a constant $C_2>0$ independent of $k$ such that on $B _\epsilon(\fkz^{(\ell)})$
\begin{equation}\label{jv_add8.2}
\mc K \leq C_2.
\end{equation}
%\begin{enumerate}
%\item[(1)] $C\inv\leq W \leq C$: this is due to the bounds on  $H_{\mathsf{f} }$
% and $W_{g_\ell}$ on $B _\epsilon(\fkz^{(\ell)})$;
%\item[(2)] $\mc K \leq C$: this follows from (1) %Theorem \ref{theorem_7.0.1}
% and Theorem \ref{theorem_3.4.2}.
%\end{enumerate}
By the convexity of $u$ and $\|u-v\|_{L^{\infty}(\Delta)}\leq \mathcal C_{1},$ we have
$$ |\partial _{2}u|\leq  \mathcal C_{1} a^{-1},\;\;\; \;\;\; \;\;\; \partial _{1} u\leq  \mathcal C_{1}a^{-1}.$$
That is
$\max_{B _\epsilon(\fkz^{(\ell)})}|z|\leq C_3.$
 Hence, by Theorem \ref{theorem_3.4.1},
 we have the regularity of $\mathsf{f} $
on $B_\epsilon(\fkz^{(\ell)})$.

\subsection{Regularity at vertices}
\label{sect_8.2}
Let $\vartheta$ be any vertex. By the  results in section \ref{sect_8.1},
  there is a bounded open set $\Omega_\vartheta\subset \CHART_\vartheta$, independent
of $k$,
such that $\vartheta \in \tau(\Omega_\vartheta)$ and the regularity of $f_\vartheta$ holds in a neighborhood of $\partial \Omega_\vartheta$.

  We omit the index $k$ if
there is  no danger of confusion.  We now quote Lemma 5.4 in \cite{CLS}.  Recall that $\phi=f_\vartheta-g_\vartheta\in C^\infty_{\mathbb T^2}(M),$
where $f_\vartheta,g_\vartheta$ are Legendre transform of $u,v$ respectively with respect to the origin $\vartheta$. Denote  $f_\vartheta$ by $f$.
Let
$$
T=\sum f^{i\bar i},\;\;\; P=\exp(\kappa W^\alpha)\sqrt{W}\Psi,
\;\;\;Q=e^{ N_1(|z|^2-A)}\sqrt{W}T.
$$

\begin{lemma}\label{lemma_8.2.1}
 Let $\Omega\subset \CHART_\vartheta$ and   $K= \mc S(f)\circ \tau_f\inv$ be the scalar curvature function
on $\t$. Suppose that on $\Omega$
\begin{equation}
\label{eqn_8.1}\max_{\tau_f(\bar \Omega)}\left(|K|+ \sum\left|\frac{\p
K}{\p \xi_i}\right|\right)\leq C_0,\;\;\;\;\;W\leq C_0,\;\;\;\;\;\;
|z|\leq C_0\end{equation}  for some constant $C_0>0$ independent of $k$. Then we may
choose
\begin{equation}\label{eqna_8.2}A=C_0^2+1,N_1=100,\alpha=\frac{1}{3},\kappa=[4C_0^\frac{1}{3}]\inv
\end{equation}
 such that
\begin{equation}\label{eqn_8.3}
\df (P+Q+ \mc C_7f ) \geq \mc C_{6}(P+Q)^2>0
\end{equation} for some positive constants $\mc C_6$ and $\mc C_7$ that
depend only on $C_0$ and $n$.
\end{lemma}

It follows from this result that $P$ and $Q$ are bounded above.   By Theorem \ref{theorem_7.0.1} we conclude that  $W$ is bounded below and above in $\Omega_\vartheta$. Then $T$ is bounded above. Therefore we have a constant $C$
such that
$$
C\inv\leq \lambda_1\leq\lambda_2\leq C,\;\;\;\;\;\;\|W\|_{C^{1}}\leq C
$$
where $\lambda_1$ and $\lambda_2$ are eigenvalues of $(f_{i\bar j})$.
By Theorem 6.7 in \cite{CLS} we have for any $U\subset\subset \Omega_\vartheta$,
$$\|f\|_{C^{6,\alpha}(U)}\leq C_{1}.$$
where $C_{1}$ is a constant depends on $\|K\|_{C^{3}(\Delta)},$ $d_E(U,\Omega_\vartheta)$ and the bound of $\Omega_\vartheta.$
 
  Then we get the interior regularity $f$ in $\Omega_\vartheta:$
\v
{\em $\phi\indexm$ uniformly $C^{6,\alpha}$-converges to
a function $\phi\in  C^{6,\alpha}_{\mathbb T^2}(M)$ with $\mathcal
S(\phi+g)=K\circ \nabla^f$.
}
\v
  We have finished the proof of
Theorem \ref{theorem_2.2.2}.
\v

By a bootstrapping argument we have $\phi\in C^{\infty}_{\mathbb T^2}(M).$ Then $\psi\in C^{\infty}(\bar \Delta).$
It follows  that $\Lambda$ is closed. This then implies Theorem
\ref{theorem_1.1}. \\


\begin{thebibliography}{s2}


\bibitem{Ab1} M. Abreu, {\em K\"ahler geometry of toric varieties and extremal metrics.} Internat. J. Math., 9(1998), 641-651.

\bibitem{ACG} V. Apostolov, D. Calderbank, P. Gauduchon,  {\em Hamiltonian 2-forms in K\"ahler geometry.} I. General theory, J. Differential Geom.,  73(2006), 359-412.

\bibitem{APS}C. Arezzo, F. Pacard, M. Singer,  {\em Extremal metrics on blow ups.} Duke Math. J., 157(2011), 1-51.

\bibitem{BA} I.J. Bakelmann, {\em Convex Analysis and Nonlinear Geometric Elliptic Equation.} Sprigner-Verlag, Berlin, 1994.

\bibitem{BU} H. Buseman,  {\em Convex Surfaces}, Interscience, New York, 1958.

\bibitem{C} E. Calabi, {\em Extremal K\"ahler metrics, Seminar on Differential Geometry,}  Ann. of Math. Studies, University Press, Series 102(1982), Princeton, 259-290.

\bibitem{C1} L.A. Caffarelli, {\em Interior $W^{2,p}$ estimates for solutions of Monge-Amp\`ere equations,}  Ann. Math., 131(1990), 135-150.

\bibitem{CG} L.A. Caffarelli, C.E. Guti\'errez, {\em Properties of the solutions of the linearized Monge-Amp\`ere equations,} Amer. J. Math., 119(1997), 423-465.




\bibitem{CLS} B.H. Chen, A.-M. Li, L.Sheng,
{\em Affine techniques on extremal metrics on toric surfaces,}  Priprint,
 ArXiv:1008.2606.

\bibitem{CLS1} B.H. Chen, A.-M. Li, L. Sheng,  {\em The Abreu equation with degenerated boundary condition,} J. Diff. Equations, 252(2012), 5235-5259.

\bibitem{CLS2} B.H. Chen, A.-M. Li, L. Sheng, {\em Interior regularization  for solutions of Abreu's equation,}  Acta Math. Sinica, 29(2013), 33-38.

\bibitem{CLS3} B.H. Chen, A.-M. Li, L. Sheng, {\em Extremal metrics on orbifold toric surfaces,} in preparation.


\bibitem{CLS4} B.H. Chen, A.-M. Li, L. Sheng, {\em On extremal Sasaki metrics on toric Sasaki manifolds,} in preparation.

\bibitem{CLS5} B.H. Chen, A.-M. Li, L. Sheng, {\em Uniform K-stability for extremal metrics on toric varieties,} J. Diff. Equations, 257(2014), 1487-1500

\bibitem{CLW} X.X. Chen, C. Lebrun, B. Weber, {\em On conformally K\"ahler, Einstein manifolds.} J. Amer. Math. Soc., 21(2008), 1137-1168.

\bibitem{CT} X.X. Chen, G. Tian, {\em Uniqueness of extremal K\"ahler metrics,}  C. R. Math. Acad. Sci. Paris, 340(2005), 287-290.

\bibitem{CT1} X.X. Chen, G. Tian, {\em Geometry of K\"ahler metrics and foliations by holomorphic discs,} Publ. Math. Inst. Hautes \'{E}tudes Sci., 107(2008), 1-107.

\bibitem{chen01} X.X. Chen,  {\em Calabi flow in Riemann surfaces revisited,}  Internat. Math. Res. Notices, 6(2001), 275-297.

\bibitem{chen05} X.X. Chen, {\em Space of K\"ahler metrics III--On the lower bound of the Calabi energy and geodesic distance,} Invent. Math., 175(2009), 453-503.

\bibitem{chen-he} X.X. Chen, W.Y. He, {\em On the Calabi flow,}  Amer. J. Math., 130(2008), 539-570.

\bibitem{chen-he-1}X.X. Chen, W.Y. He, {\em The Calabi flow on K\"ahler surface with bounded Sobolev constant (I),} Math. Ann., 354(2012), 227-261.

\bibitem{cw} X.X. Chen, B. Weber, {\em Moduli spaces of critical Riemannian metrics with $L^{n\over 2}$ norm curvature bounds,}  Adv. Math., 226(2011), 1307-1330.

\bibitem{De} T. Delzant, {\em Hamiltioniens periodiques et image convex de l'application moment,} Bull. Soc. Math. France, 116(1988), 315-339.

\bibitem{D-Z-Z} S. Dinew, X. Zhang, X.W. Zhang, {\em The $C^{2,\alpha}$  estimate of complex Monge-Amp\`ere equation.} Indiana Univ. Math. J., 60(2011), 1713-1722.

\bibitem{D1} S.K. Donaldson, {\em Scalar curvature and stability of toric varietties.} J. Differential Geom., 62(2002), 289-349.

\bibitem{D2} S.K. Donaldson, {\em Constant scalar curvature metrics on toric surfaces.}  Geom. Funct. Anal.,  19(2009), 83-136.

\bibitem{D3} S.K. Donaldson, {\em Interior estimates for solutions of Abreu's equation.} Collectanea Math., 56(2005), 103-142.

\bibitem{D4} S.K. Donaldson, {\em Extremal metrics on toric surfaces: a continuity method,}  J. Differential Geom., 79(2008), 389-432.

\bibitem{Fu} W. Fulton, {\em Introduction to Toric Varieties,} Ann. of Math. Studies, University Press, Series 131(1993), Princeton.

\bibitem{G} C.E. Guti\'errez, {\em The Monge-Amp\`ere Equations,} Progress in Nonlinear Differential Equations and Their Applications, Vol.44, Birkhauser, 2001.

\bibitem{Gu} V. Guillemin, {\em Moment Maps and Combinatorial Invariants of Hamiltonian $T^n$-Spaces,}  Progress in Mathematics, Boston, Series 122, 1994.

\bibitem{Gu1} B. Guan, {\em The Dirichlet problem for complex Monge-Amp\`ere equations and applications,}  in ``Trends in Partial Differential Equations", ALM10, 53-97, Higher Education Press and IP, Beijing-Boston, 2009.

\bibitem{Gu2} B. Guan, {\em The Dirichlet problem for complex Monge-Amp\`ere equations and regularity of the pluri-complex Green function,} Comm. Anal. Geom., 6(1998), 687-703; a correction, 8(2000), 213-218.

\bibitem{H}   Q. Han, {\em A Course in Nonlinear Elliptic Differential Equations}, Lecture Notes, Peking University, 2013.

\bibitem{H1} W.Y. He, {\em Remarks on the existence of bilaterally symmetric extremal K\"ahler metrics on $\mathbb{CP}^2\sharp 2\bar{\mathbb{CP}^2},$}  Internat. Math. Res. Notices, (2007), Art. ID rnm127, 13pp.

\bibitem{L-J} F. Jia, A.-M. Li, {\em Complete K\"ahler affine manifolds,} Preprint,  ArXiv:1008.2604.

\bibitem{L-1} A.-M. Li, {\em Calabi conjecture on hyperbolic affine hyperspheres (2),} Math. Ann., 293(1992), 485-493.

\bibitem{L-J-1} A.-M. Li, F. Jia, {\em A Bernstein property of affine maximal hypersurfaces,}  Ann. Global Anal. Geom., 23(2003), 359-372.

\bibitem{L-J-2} A.-M. Li, F. Jia, {\em Euclidean complete affine surfaces with constant affine  mean curvature,} Ann. Global Anal. Geom., 23(2003), 283-304.

\bibitem{L-J-3} A.-M. Li, F. Jia, {\em The Calabi conjecture
on affine maximal surfaces}, Result. Math., 40(2001), 265-272.

\bibitem{L-J-4} A.-M. Li, F. Jia, {\em A Bernstein properties of some fourth order partial differential equations}, Result. Math., 56(2009), 109-139.

\bibitem{LeBrun-Simanca} C. LeBrun, S. R. Simanca,  {\em Extremal K\"ahler metrics and complex deformation theory,} Geom. Funct. Anal., 4(1994), 298-336.

\bibitem{L-S-Z} A.-M. Li, U. Simon, G. Zhao, {\em Global Affine Differential Geometry of Hypersurfaces,} Walter de Gruyter, Berlin, New York, 1993.

\bibitem{L-X1} A.-M. Li,  R.W. Xu,  {\em A rigidity theorem for affine K\"{a}hler-Ricci flat graph,}  Result. Math., 56(2009), 141-164.

\bibitem{L-X2} A.-M. Li,  R.W. Xu, {\em A cubic form differential inequality  with applications to affine K\"{a}hler  Ricci flat  manifolds}, Result. Math., 54(2009), 329-340.

\bibitem{L-X-S-J} A.-M. Li,  R.W. Xu,  U. Simon, F. Jia, {\em Affine Bernstein Problems and Monge-Amp\`ere Equations,}  World Scientific, 2010.

\bibitem{L-Ss} C. LeBrun, S.R. Simanca, {\em Extremal K\"ahler
metrics and complex deformation theory,} Geom. Funct. Anal.,
4(1994), no.3, 298-336.

\bibitem{M} T. Mabuchi, {\em Extremal metrics and stabilities on polarized manifolds,}  International Congress of Mathematicians. Vol.II, 813-826, Eur. Math. Soc., Z\"urich, 2006.

\bibitem{P} A.V. Pogorelov, {\em The Minkowski Multidimensional Problem,}  John Wiley, 1978.

\bibitem{L-S} L. Sheng, G. Chen, {\em A Bernstein property of some affine K\"ahler scalar flat graph,} Result. Math., 56 (2009), 165-175.

\bibitem{Sz} G. Sz\'{e}kelyhidi, Extremal metrics and K-stability (PhD
thesis).

\bibitem{T} G. Tian, {\em Extremal metrics and geometric stability,} Houston J. Math., 28(2002), 411-432.

\bibitem{T1} G. Tian, {\em K\"ahler-Einstein metrics with postive scalar curvature,} Invent. Math., 130(1997), 1-39.

\bibitem{T2} G. Tian, {\em Canonical Metrics in Kahler Geometry.}  Lectures in Mathematics Eth Zurich,  Birkh\"auser, Basel, 2000.

 \bibitem{T3} G. Tian, {\em On Calabi's conjecture for complex surfaces with positive first Chern class,} Invent. Math., 101(1990), 101-172.

 \bibitem{T4} G. Tian, {\em The $K$-energy on hypersurfaces and stability,} Comm. Anal. Geom.,
 2(1994), 239-265.

  \bibitem{T5} G. Tian, {\em Bott-Chern forms and geometric stability,}
Discrete Contin. Dynam. Systems, 6(2000), 211-220.

\bibitem{TW} N.S. Trudinger, X.J. Wang, {\em The affine Plateau problem,}  J. Amer. Math. Soc., 18(2005), 253-289.

\bibitem{Yau} S.T. Yau, {\em On the Ricci curvature of compact K\"ahler manifold and the complex Monge-Amp\`ere equation. I,}  Comm. Pure Appl. Math., 31(1978), 339-411.


\bibitem{ZZ1}B. Zhou, X.H. Zhu,  {\em Relative $K$-stability and modified $K$-energy on toric manifolds,}  Adv. Math., 219(2008), 1327-1362.

\bibitem{ZZ2} B. Zhou, X.H. Zhu, {\em Minimizing weak solutions for Calabi's extremal metrics on toric manifolds,} Calc. Var. \& PDE., 32(2008), 191-217.

\end{thebibliography}
\end{document}